\documentclass[12pt,oneside,notitlepage]{amsart}
\usepackage{graphicx}
\usepackage{amssymb}
\usepackage{xcolor}
\newtheorem{thm}{Theorem}[section]
\newtheorem{cor}[thm]{Corollary}
\newtheorem{lem}[thm]{Lemma}
\newtheorem{prop}[thm]{Proposition}
\theoremstyle{definition}

\theoremstyle{remark}
\newtheorem{rem}[thm]{Remark}
\newtheorem{ex}[thm]{Example}
\numberwithin{equation}{section}
\newcommand{\norm}[1]{\left\Vert#1\right\Vert}
\newcommand{\abs}[1]{\left\vert#1\right\vert}
\newcommand{\set}[1]{\left\{#1\right\}}

\newcommand{\eps}{\varepsilon}
\newcommand{\To}{\longrightarrow}
\newcommand{\Hilb}{\textrm{Hilb}}
\newcommand{\FS}{FS}
\newcommand{\Aut}{\textrm{Aut}}
\newcommand{\Herm}{\textrm{Herm}}
\newcommand{\tr}{\textrm{tr}}
\newcommand{\Sym}{\textrm{Sym}}
\newcommand{\dist}{\textrm{dist}}
\newcommand{\Exp}{\textrm{Exp}}
\newcommand{\Lie}{\textrm{Lie}}
\newcommand{\CP}{\mathbb{CP}}
\newcommand{\ddbar}{\sqrt{-1}\,\partial\bar{\partial}}

\newcommand{\R}{\mathbb{R}}
\newcommand{\C}{\mathbb{C}}
\newcommand{\M}{\mathcal{M}}
\newcommand{\K}{\mathcal{K}}
\newcommand{\B}{\mathcal{B}}
\newcommand{\Nor}{\mathcal{N}}
\newcommand{\Zt}{\widetilde{Z}}
\newcommand{\pp}{\mathfrak{p}}

\begin{document}

\title[Numerical algorithm for balanced metrics]{Numerical algorithm for finding balanced metrics on polarized manifolds}%

\author[R. Seyyedali]{Reza Seyyedali}
\address {School of Mathematics, Institute for Research in Fundamental Sciences (IPM), P.O. Box: 19395-5746, Tehran, Iran}
\email{rseyyedali@ipm.ir}
\address{Department of Mathematics, New Uzbekistan University, Tashkent,  Uzbekistan}
\email{r.seyyedali@newuu.uz}

\thanks{}%

\begin{abstract}
In \cite{D3}, Donaldson defines a dynamical system $T$ on the space of Fubini--Study metrics of a polarized manifold $(X,L)$, whose fixed points are the balanced metrics. In \cite{S}, Sano proves that the iterates of $T$ converge to a balanced metric, provided that a balanced metric exists and the automorphism group of $(X,L)$ is discrete modulo scalars. We remove the assumption on the automorphism group. Moreover, we show that the convergence is exponential with the optimal rate $\rho_*$ given by the linearization of $T$ at a balanced metric, which extends a result of Ioos and Polterovich (\cite{IP}) for discrete automorphism groups. We compute $\rho_*$ for $(\CP^n,\mathcal{O}(d))$ and $(\CP^1\times\CP^1,\mathcal{O}(d,d))$; for $\CP^1$ this confirms the rate obtained numerically by Sherman (\cite{Sh}).
\end{abstract}
\maketitle
\pagestyle{plain} \pagenumbering{arabic}

\section{Introduction}

A central notion in geometric invariant theory is the concept of stability. This notion is crucial to form quotient spaces in algebraic geometry. An important example of such quotients is the moduli space of projective algebraic varieties. The most well known constructions are through Hilbert points and Chow points of projective varieties, and the corresponding stability conditions are Hilbert stability and Chow stability. By a general theme that relates GIT quotients to symplectic quotients, one expects that these notions of stability correspond to canonical metrics on the stable objects. In \cite{Zh}, Zhang introduced the concept of balanced embedding and proved that the existence of a balanced embedding of a polarized variety is equivalent to the stability of its Chow point. Zhang's result has been reproven by Luo in \cite{L} and by Phong and Sturm in \cite{PS}. In \cite{D3}, Donaldson introduces a dynamical system on the space of Fubini--Study metrics in order to approximate balanced metrics numerically. In \cite{S}, Sano proves that this dynamical system converges to the balanced metric, provided that a balanced metric exists and the automorphism group is discrete modulo scalars.

To state the precise result, let $X$ be a compact connected complex manifold of dimension $n$ and $L$ be a very ample line bundle on $X$. Let $W=H^{0}(X,L)$ and $N=\dim W$. We denote the space of Hermitian inner products on $W$ by $\M$ and the space of Hermitian metrics $h$ on $L$ with positive curvature form $\omega_h=-\ddbar\log h$ by $\K$. For $h\in\K$, we put $d\mu_h=\omega_h^n/n!$ and $V=\int_X d\mu_h$. Note that $V$ does not depend on $h$ and only depends on the polarization $L$. Following Donaldson (\cite{D3}), we define the maps $\Hilb:\K\to\M$ and $\FS:\M\to\K$ as follows.
\begin{itemize}
\item For $h\in\K$, the inner product $\Hilb(h)$ on $W$ is defined by
$$\Hilb(h)(s,t)=\frac{N}{V}\int_X h(s,t)\,d\mu_h.$$
\item For $H\in\M$, the metric $\FS(H)$ is the unique metric on $L$ such that
$$\sum_{i=1}^N|s_{i}|_{\FS(H)}^2=1,$$
where $s_{1},\dots,s_{N}$ is an orthonormal basis of $W$ with respect to $H$.
\end{itemize}
Let $T=\Hilb\circ\FS:\M\to\M$. An inner product $H\in\M$ is called balanced if $T(H)=H$. For $h\in\K$, the Bergman kernel of $h$ is the function $\rho_h=\sum|t_i|_h^2$, where $t_1,\dots,t_N$ is an orthonormal basis of $W$ with respect to $\Hilb(h)$. It is easy to see that $H$ is balanced if and only if the Bergman kernel of the Hermitian metric $\FS(H)$ on $L$ is constant. We denote the set of balanced inner products by $\B$. By the work of Zhang and Luo (\cite{Zh}, \cite{L}), $\B$ is nonempty if and only if the Chow point of the embedding of $X$ by sections of $L$ is polystable. In \cite{D1}, Donaldson shows that if $\Aut(X,L)/\C^*$ is discrete, then balanced metrics are unique up to scaling. In general, $\B$ is an orbit of the group of automorphisms of $(X,L)$ (Proposition \ref{propB}). The space $\M$ can be naturally identified with the symmetric space $GL(N,\C)/U(N)$, and therefore it carries a natural Riemannian metric with corresponding distance 
$$\dist(H,H')=\Big(\sum_{i=1}^{N}(\log\mu_{i})^{2}\Big)^{\frac{1}{2}},$$
where $\mu_{1},\dots,\mu_{N}$ are the eigenvalues of $H^{-1}H'$. Our main theorem is the following

\begin{thm}\label{mainthm}
Suppose that $(X,L)$ admits a balanced metric. Then for any $H_{0}\in \M$, the sequence $T^m(H_{0})$ converges to a balanced inner product $H_{\infty}\in\B$. Moreover, there exists $\rho_*\in(0,1)$ depending only on $(X,L)$ such that the following holds. For any $R>0$, there is a constant $C=C(R)$ such that
$$\dist\big(T^m(H_0),H_{\infty}\big)\leq C\rho_*^m\qquad\textrm{for all } m\geq 0,$$
whenever $\dist(H_0,\B)\leq R$.
If $T$ is not the identity map, then the rate $\rho_*$ is optimal.
\end{thm}

Here $\rho_*$ is the largest eigenvalue of the linearization of $T$ at a balanced metric, restricted to the directions normal to $\B$ (see Theorem \ref{thmlin}).  Note that $H_\infty$ depends on $H_0$. So, in the case of $\B\neq \M$, the constant $C$ cannot be chosen independent of $R$ (Remark \ref{remC}), and the rate $\rho_*$ cannot be replaced by any $\lambda<\rho_*$ (Proposition \ref{propsharp}). It is worth noting that  $T$ is  the identity map, i.e. that $\B=\M$, only in the case of  $(\CP^n,\mathcal{O}(1))$.

By the celebrated theorem of Donaldson (\cite{D1}), if a polarization $L$ admits a constant scalar curvature K\"ahler metric and $\Aut(X,L)/\C^*$ is discrete, then $L^k$ carries balanced metrics, unique up to scaling, for $k\gg0$. In \cite[Theorem 1.2(2)]{IP}, Ioos and Polterovich prove that, under these assumptions and for $k\gg0$, the iterates of $T$ for the polarization $L^k$ converge exponentially to a balanced metric. By the proof of \cite[Theorem 1.2]{IP}, the rate is the largest eigenvalue of the linearization of $T$ different from $1$, i.e. the number $\rho_*$ of Theorem \ref{mainthm}, and they compute its asymptotics as $k\to\infty$ (see Section \ref{secasym}). Note that Theorem \ref{mainthm} assumes neither condition, and that the constant $C$ depends only on $R$.  The optimality of the rate (Proposition \ref{propsharp}) is expected from the toy model of \cite[Section 2.1.1]{D3}. We now compare Theorem \ref{mainthm} with other related results. In \cite[Theorem 1.2]{Io1}, Ioos proves the exponential convergence of the anticanonical analogue of $T$ on Fano manifolds with discrete automorphism group, and computes the asymptotics of the optimal rate. For the map $T_\nu$ of \cite{D3}, in which a measure $\nu$ on $X$ is fixed, Ioos, Kaminker, Polterovich and Shmoish prove exponential convergence under conditions on $\nu$ which do not involve automorphisms (\cite[Theorem 4.4]{IKPS}). In \cite[Theorem 4.14]{Be}, Berman proves the convergence of the canonical and anticanonical analogues of $T$. In the anticanonical case, he assumes that $X$ has no nonzero holomorphic vector fields and that an anticanonically balanced metric exists.

In \cite{S}, Sano shows that $\Zt(H_m)$ is decreasing, where $\Zt$ is a functional on $\M$ introduced by Donaldson whose critical points are exactly the balanced metrics. Then he uses the properness of $\Zt$ to show that the sequence $H_m=T^m(H_0)$ is bounded. Note that Sano does not use any assumption on the automorphism group in order to prove that $\Zt(H_m)$ is decreasing. However, if $\Aut(X,L)/\C^*$ is not discrete, then $\Zt$ is constant on $\B$, which is noncompact even modulo scaling. Therefore $\Zt$ is neither strictly convex nor proper, even modulo scaling, and the boundedness of $H_m$ is not clear. In order to overcome this issue, we first show that $\B$ is a closed totally geodesic submanifold of $\M$ (Proposition \ref{propgeomB}). This gives a substitute for the properness of $\Zt$, namely properness modulo the automorphism group, which follows from convexity alone (Lemma \ref{lemproper}). Note that such properness can also be proven using the Kempf--Ness theorem (Remark \ref{remKN}). Using this, we prove that the sequence $H_m$ approaches $\B$ (Proposition \ref{propconv}). It remains to rule out the situation in which $H_m$ goes to infinity along $\B$. We give two essentially different proofs of this. The first proof uses the convexity of $\Zt$ to prove a {\L}ojasiewicz type inequality near $\B$ (Lemma \ref{lemloj}). The second proof relies on the linearization of $T$ at a balanced metric $H$. We show that
$$d_HT=I-\frac{N}{V}\nabla^2\Zt(H),$$
where $\nabla^2\Zt$ is the Hessian of $\Zt$. Hence $d_HT$ is self-adjoint with eigenvalues in $(0,1]$, and the eigenvalue $1$ occurs exactly in the directions tangent to $\B$ (Theorem \ref{thmlin}).

Combining Theorem \ref{thmlin} with the fact that the automorphism group acts on $\M$ by isometries, we show that near $\B$ the distance to $\B$ contracts at each step, while the displacement along $\B$ is controlled by the distance to $\B$. The second proof gives the rate $\rho_*$. Using the same reduction to the normal directions of $\B$, we show that this rate cannot be improved (Proposition \ref{propsharp}). In \cite{IP}, Ioos and Polterovich also compute the linearization of $T$ at a balanced metric, in terms of the Berezin transform (see Remark \ref{remberezin}).

Both proofs follow well-known methods in analsis. The first proof is the {\L}ojasiewicz type argument for descent methods (c.f. \cite{AMA}). It combines the sufficient decrease of Lemma \ref{lemdescent} with a gradient inequality with exponent $\frac12$ (Lemma \ref{lemloj}), which holds since the Hessian of $\Zt$ is nondegenerate in the directions normal to $\B$. The second proof is the analysis of a smooth map near a manifold of fixed points, whose linearization is a contraction in the normal directions (c.f. \cite{Au}). The new point is that $\B$ is noncompact. The uniformity of the estimates along $\B$ comes from the action of $G$ on $\M$ by isometries.

In \cite{Sh}, Sherman studies numerically the convergence of $T$-iterations on the Riemann sphere, and observes that for $S^1$-invariant metrics on $(\CP^1,\mathcal{O}(k))$ the rate of convergence is $\frac{(k-1)(k+6)}{(k+2)(k+3)}$ (\cite[(9)]{Sh}). Note that the automorphism group is continuous in that case, which is covered by Theorem \ref{mainthm}, without the assumption of $S^1$-invariance. In Section \ref{secexample}, we compute $\rho_*$ explicitly for $(\CP^n,\mathcal{O}(d))$ and for $(\CP^1\times\CP^1,\mathcal{O}(d,d))$. For $n=1$, this gives Sherman's formula (Remark \ref{remsherman}). For $\CP^1\times\CP^1$, we get $\rho_*=1-\frac{4}{(d+2)^2}$. Note that in this case the slowest direction is a mixed mode, which has no counterpart on either factor (Proposition \ref{propproduct}).

In \cite{Se1}, we proved the analogue of Sano's theorem for balanced metrics on simple vector bundles in the sense of Wang (\cite{W}). In \cite{Se3}, we remove the simplicity assumption by introducing the Thompson metric on $\M$, given by $d_{Th}(H,H')=\max_i\abs{\log\mu_i}$, where $\mu_1,\dots,\mu_N$ are the eigenvalues of $H^{-1}H'$. The map $T$ for vector bundles preserves the partial order on $\M$ and therefore it is distance non-increasing. In the manifold case, this approach fails. Indeed, the map $T$ is neither order preserving nor distance non-increasing with respect to the Thompson metric, even near a balanced metric (Remark \ref{remthompson}).

The outline of the paper is as follows: In Section $2$, we go over the functional $\Zt$ and the geometry of $\M$, and we prove that $\B$ is a closed totally geodesic submanifold of $\M$. In Section $3$, we prove that the sequence $T^m(H_0)$ approaches $\B$. Section $4$ is devoted to two proofs of Theorem \ref{mainthm} and to the optimality of the rate $\rho_*$. In Section $5$, we compute $\rho_*$ for projective spaces and for $\CP^1\times\CP^1$, and we discuss the behaviour of $\rho_*$ for high powers of the polarization.

\medskip
\noindent\textbf{Acknowledgements:} I used the AI assistant Claude (Opus 5.5, Anthropic) in the preparation of this paper. Claude helped with the calculations in Example \ref{excp1} and in the examples of Section \ref{secexample}, and it was used to expedite the proofs of Proposition \ref{propgeomB} and of the results in Appendix \ref{appfutaki}. Claude was also used to review a draft of the paper and to check numerically some cases of the rates of convergence computed in the examples of Section \ref{secexample}, and it suggested the references \cite{AMA},and \cite{Au} . It also assisted in identifying and correcting minor errors and in filling minor gaps. I emphasize that the main results and their proofs, the overall structure of the paper and the vast majority of the examples are my own work, with some help from Claude in expediting calculations. I have checked all of the mathematics and the references, and I take full responsibility for the content of the paper.

\section{Balanced metrics and the functional $\widetilde{Z}$}

Throughout the paper, $X$, $L$, $W$, $N$, $\M$, $\K$, $d\mu_h$ and $V$ are as in the introduction. Note that for $h\in\K$ and a smooth real function $\varphi$ on $X$, we have $\omega_{e^{-\varphi}h}=\omega_h+\ddbar\varphi$.

Let $H\in \M$ and $\underline{s}=(s_{1},\dots,s_{N})$ be an $H$-orthonormal basis of $W$. Since $L$ is very ample, the basis $\underline{s}$ defines an embedding
$$\iota_{\underline{s}}=[s_{1}:\dots:s_{N}]:X\To \CP^{N-1}$$
satisfying $\iota_{\underline{s}}^{*}\mathcal{O}(1)=L$. The metric $\FS(H)$ is the unique Hermitian metric on $L$ such that
\begin{equation}\label{eqFS}
\sum_{i=1}^{N}|s_{i}|^{2}_{\FS(H)}=1.
\end{equation}
Note that a change of orthonormal basis corresponds to an element of $U(N)$ and therefore $\FS(H)$ does not depend on the choice of the orthonormal basis. Moreover, $\FS(H)$ is the pull back of the Fubini--Study metric of $\mathcal{O}(1)$ under $\iota_{\underline{s}}$. Hence $\FS(H)\in \K$ and $\omega_{\FS(H)}=\iota_{\underline{s}}^{*}\omega_{FS}$, where $\omega_{FS}=\ddbar\log|z|^2$ is the Fubini--Study form on $\CP^{N-1}$. We write $\iota_H$ for $\iota_{\underline{s}}$ when the choice of basis is irrelevant. If $H$ is balanced, then $\iota_{H}(X)$ is a balanced embedding in the sense of \cite{Zh} and \cite{L}. For any real number $\alpha$, we have the scaling identities
$$\FS(e^{\alpha}H)=e^{\alpha}\FS(H), \quad \Hilb(e^{\alpha}h)=e^{\alpha}\Hilb(h),\quad T(e^{\alpha}H)=e^{\alpha}T(H).$$

For $H\in \M$, we denote by $H^{-1}T(H)$ the endomorphism of $W$ defined by $T(H)(s,t)=H(H^{-1}T(H)s,t)$, and we define
\begin{equation}\label{eqP}
P(H)=H^{-1}T(H)-I.
\end{equation}

\begin{lem}\label{lemtrace}
For any $H\in\M$, the endomorphism $H^{-1}T(H)$ is self-adjoint and positive with respect to $H$, and $\tr\big(H^{-1}T(H)\big)=N$. In particular, the eigenvalues $p_1,\dots,p_N$ of $P(H)$ satisfy $-1<p_i<N-1$ and $\sum p_i=0$.
\end{lem}

\begin{proof}
Let $s_1,\dots,s_N$ be an $H$-orthonormal basis and $h=\FS(H)$. Since $T(H)$ is an inner product, $H^{-1}T(H)$ is self-adjoint and positive with respect to $H$. By \eqref{eqFS}, we have
$$\tr\big(H^{-1}T(H)\big)=\sum_{i=1}^N T(H)(s_i,s_i)=\frac{N}{V}\int_X\sum_{i=1}^N|s_i|_h^2\,d\mu_h=N.$$
Hence the eigenvalues $1+p_i$ of $H^{-1}T(H)$ are positive and sum to $N$.
\end{proof}
Following Aubin and Yau, we define the functional $I:\K \to \R$ by its variation along a smooth path $h_{t}$ in $\K$ as follows:
\begin{equation}\label{eqI}
\frac{d}{dt}I(h_{t})=\int_{X}\frac{d}{dt}\log\Big(\frac{h_{t}}{h_{0}}\Big)\,d\mu_{h_{t}}.
\end{equation}
Note that the functional $I$ is defined up to an additive constant. Fixing a nonzero element of $\bigwedge^{N}W$, we define $\log\det:\M\To\R$. Following Donaldson, we define
$$Z=-I\circ \FS:\M\To\R, \qquad \Zt=Z+\frac{V}{N}\log\det.$$
The scaling identities imply that $\Zt(e^{\alpha}H)=\Zt(H)$.

Let $H\in \M$ and $a$ be an $H$-self-adjoint endomorphism of $W$. Let $s_{1},\dots,s_{N}$ be an $H$-orthonormal basis such that $as_{i}=\lambda_{i}s_{i}$. We define the path $H_{t}\in \M$ by
\begin{equation}\label{eqgeod}
H_{t}(s_{i},s_{j})=\delta_{ij}e^{t\lambda_{i}}.
\end{equation}
These are the geodesics of the symmetric space $\M$ (see below).

\begin{lem}\label{lemfirst}
With the above notation,
$$\frac{d}{dt}\Big|_{t=0}\Zt(H_{t})=-\frac{V}{N}\tr\big(aP(H)\big).$$
In particular, $H$ is a critical point of $\Zt$ if and only if $H$ is balanced.
\end{lem}

\begin{proof}
Put $h=\FS(H)$ and $h_{t}=\FS(H_{t})$. Since $e^{-\frac{t\lambda_{1}}{2}}s_{1},\dots,e^{-\frac{t\lambda_{N}}{2}}s_{N}$ is $H_{t}$-orthonormal, we have $h_{t}=e^{-\varphi_{t}}h$, where
\begin{equation}\label{eqphi}
e^{\varphi_{t}}=\sum_{i=1}^{N}e^{-t\lambda_{i}}|s_{i}|_{h}^{2}.
\end{equation}
Hence $\varphi_{0}=0$ and $\dot{\varphi}_{0}=-\sum\lambda_{i}|s_{i}|_{h}^{2}$. By \eqref{eqI}, we have,
\begin{align*}
\frac{d}{dt}\Big|_{t=0}\Zt(H_{t})&=\int_{X}\dot{\varphi}_{0}\,d\mu_{h}+\frac{V}{N}\sum \lambda_i=-\sum\lambda_{i}\int_{X}|s_{i}|_{h}^{2}d\mu_{h}+\frac{V}{N}\sum \lambda_i\\
&=-\frac{V}{N}\sum\lambda_{i}\,\Big(T(H)(s_{i},s_{i})-1\Big)=-\frac{V}{N}\tr\big(aP(H)\big).
\end{align*}

Hence $H$ is a critical point of $\Zt$ if and only if $P(H)=0$ i.e. $H$ is balanced.
\end{proof}

The following lemma is due to Sano (\cite [Section $2$]{S}).

\begin{lem}\label{lemmono}
For any $H\in \M$,
\begin{enumerate}
    \item $Z(T(H))\leq Z(H),$
    \item $\log\det T(H)\leq \log\det H,$
    \item $\Zt(T(H))\leq \Zt(H).$
\end{enumerate}

\end{lem}

The proof of the second part of Lemma \ref{lemmono} uses the arithmetic-geometric mean inequality:
$$\det(H^{-1}T(H))\leq \Big(\frac{\tr(H^{-1}T(H))}{N}\Big)^{N}=1.$$
Therefore equality  holds if and only if $T(H)=H$. A simple calculation gives the following quantitative version of Lemma \ref{lemmono}.

\begin{lem}\label{lemdescent}
For any $H\in \M$, we have
\begin{align*}
\log\det H-\log\det T(H)&\geq\frac{1}{2N}\norm{P(H)}^{2},\\
\Zt(H)-\Zt(T(H))&\geq\frac{V}{2N^{2}}\norm{P(H)}^{2}.
\end{align*}
Moreover, if $\norm{P(H)}_{op} \leq \frac{1}{2}$ (say), then 
$$\dist(H,T(H)) \leq \sqrt2 \norm{P(H)}.$$
\end{lem}

\begin{proof}
Let $p_{1},\dots,p_{N}$ be the eigenvalues of $P(H)$. Lemma \ref{lemtrace} implies that $p_{i} \in (-1,N-1)$ and $\sum p_{i}=0$. Hence
$$\log\det H-\log\det T(H)=-\log\det(I+P(H))=\sum_{i=1}^{N}\big(p_{i}-\log(1+p_{i})\big).$$
  For $-1 <x\leq N-1,$ we have $$x-\log (1+x) \geq \frac{x^2}{2N}.$$  This proves the first inequality. The second one follows from the first inequality and Lemma \ref{lemmono}(1), since $\Zt=Z+\frac{V}{N}\log\det$.
  For the last statement, note that $\log^2(1+x)\leq 2x^2$ for $\abs{x}\leq\frac12$. Hence
  $$\dist(H,T(H))^2=\sum \log^2(1+p_i)\leq 2\sum |p_i|^2=2\norm{P(H)}^2.$$
\end{proof}

We now compute the second variation of $Z$ along the paths \eqref{eqgeod}. This is the second variation formula of Zhang (\cite{Zh}, see also \cite{D2}, \cite{PS}). For the reader's convenience, we include a proof. Let $H$, $a$, $\lambda_{i}$, $s_{i}$ and $H_{t}$ be as in \eqref{eqgeod}. Using the basis $s_{1},\dots,s_{N}$, we identify $\mathbb{P}(W^{*})$ with $\CP^{N-1}$ and define the function
$$\mu_{a}([z])=\frac{\sum\lambda_{i}|z_{i}|^{2}}{\sum|z_{i}|^{2}}$$
on $\CP^{N-1}$. More generally, for a Hermitian matrix $b$ we put $\mu_b([z])=z^*bz/|z|^2$, so that $\mu_{a^2}([z])=\sum\lambda_i^2|z_i|^2/\sum|z_i|^2$. Let $X_{t}=\iota_{H_{t}}(X)\subset\CP^{N-1}$ be the image of $X$ under the embedding given by the $H_{t}$-orthonormal basis $e^{-\frac{t\lambda_{i}}{2}}s_{i}$. Then $X_{t}=e^{-\frac{ta}{2}}(X_{0})$. We denote by $\omega_{FS}=\ddbar\log|z|^2$ the Fubini--Study form on $\CP^{N-1}$, by $dvol_{FS}=\omega_{FS}^n/n!$ the induced volume form on $X_t$, by $\xi_{a}=\nabla^{1,0}\mu_a$ the $(1,0)$-gradient of $\mu_{a}$ with respect to $\omega_{FS}$, and by $\xi_{a}^{\perp}$ its component normal to $X_{t}$. Note that $\xi_a$ is the holomorphic vector field on $\CP^{N-1}$ induced by $a$, i.e. it generates the action of $e^{\zeta a}$, $\zeta\in\C$.

\begin{lem}\label{lemsecond}( \cite{Zh}, see also \cite{D2}, \cite{PS})
For any $t\in\R$, we have
$$\frac{d^{2}}{dt^{2}}Z(H_{t})=\int_{X_{t}}\big|\xi_{a}^{\perp}\big|^{2}_{\omega_{FS}}\,dvol_{FS}.$$
In particular $Z$ and $\Zt$ are convex along the paths \eqref{eqgeod}. Moreover $\frac{d^{2}}{dt^{2}}\big|_{t=0}Z(H_{t})=0$ if and only if $e^{ta}$ preserves $X_{0}$ for all $t\in\R$.
\end{lem}

\begin{proof}
 It suffices to consider $t=0$. Indeed, $s\mapsto H_{t+s}$ is the path \eqref{eqgeod} starting at $H_t$, defined by the same endomorphism $a$ and the $H_t$-orthonormal basis $e^{-\frac{t\lambda_1}{2}}s_1,\dots,e^{-\frac{t\lambda_N}{2}}s_N$. Hence the case $t=0$, applied at $H_t$, gives the general case. Let $h=\FS(H)$. With the notation of the proof of Lemma \ref{lemfirst}, \eqref{eqphi} gives
$$\dot{\varphi}_{t}=-\frac{\sum\lambda_{i}e^{-t\lambda_{i}}|s_{i}|_{h}^{2}}{\sum e^{-t\lambda_{j}}|s_{j}|_{h}^{2}}=-\mu_{a}\circ\iota_{H_{t}},\qquad \ddot{\varphi}_{0}=\big(\mu_{a^2}-\mu_a^2\big)\circ\iota_{H}.$$
By \eqref{eqI}, we have $\frac{d}{dt}Z(H_{t})=\int_{X}\dot{\varphi}_{t}\,d\mu_{h_{t}}$, where $h_t=\FS(H_t)=e^{-\varphi_t}h$. Since $\omega_{h_t}=\omega_h+t\ddbar\dot\varphi_0+O(t^2)$, we have $\frac{d}{dt}\big|_{t=0}d\mu_{h_{t}}=\Delta_h\dot\varphi_0\,d\mu_h$, where $\Delta_hf\,d\mu_h=\ddbar f\wedge\omega_h^{n-1}/(n-1)!$. Integrating by parts, we get
$$\frac{d^{2}}{dt^{2}}\Big|_{t=0}Z(H_{t})=\int_X\big(\mu_{a^2}-\mu_a^2\big)\circ\iota_H\,d\mu_h-\int_X|\partial\hat a|^2_{\omega_h}\,d\mu_h,$$
where $\hat a=\mu_a\circ\iota_H$ and $|\partial f|^2_{\omega_h}d\mu_h=\sqrt{-1}\,\partial f\wedge\bar\partial f\wedge\omega_h^{n-1}/(n-1)!$ for real functions $f$.

Let $b$ be a Hermitian matrix. In the affine coordinates $w_j=z_j/z_0$ around $[1:0:\dots:0]$, we have $\omega_{FS}=\sqrt{-1}\sum dw_j\wedge d\bar w_j$ at $w=0$, and a direct computation gives
$$\partial\mu_b=\sum_{j\geq1}b_{0j}\,dw_j,\qquad \mu_{b^2}-\mu_b^2=\sum_{j\geq1}|b_{0j}|^2$$
at $w=0$. Since $U(N)$ acts transitively on $\CP^{N-1}$ by isometries and $\mu_{u^{-1}bu}=\mu_b\circ u$ for $u\in U(N)$, we conclude that $\mu_{a^2}-\mu_a^2=|\partial\mu_a|^2_{\omega_{FS}}=|\xi_a|^2_{\omega_{FS}}$ on $\CP^{N-1}$. On the other hand, $\omega_h=\iota_H^*\omega_{FS}$, so the $(1,0)$-gradient of $\hat a$ with respect to $\omega_h$ is the component of $\xi_a$ tangent to $X_0$. Therefore
$$\frac{d^{2}}{dt^{2}}\Big|_{t=0}Z(H_{t})=\int_{X_{0}}\Big(|\xi_a|^2_{\omega_{FS}}-|\xi_a^{T}|^2_{\omega_{FS}}\Big)\,dvol_{FS}=\int_{X_{0}}\big|\xi_{a}^{\perp}\big|^{2}_{\omega_{FS}}\,dvol_{FS},$$
where $\xi_a^T$ is the component of $\xi_a$ tangent to $X_0$. 
Finally, note that the second derivative vanishes at $t=0$ if and only if $\xi_{a}$ is tangent to $X_{0}$. This implies that the flow of $\xi_a$ preserves $X_0.$ Conversely, if $e^{ta}$ preserves $X_0$ for all $t\in\R$, then $\textrm{Re}\,\xi_a$ is tangent to $X_0$. Since $X_0$ is a complex submanifold, $\xi_a$ is tangent to $X_0$, and therefore $\xi_a^\perp=0$.
\end{proof}

\begin{cor}\label{cormin}
Suppose that $H_{b}\in \M$ is balanced. Then $\Zt(H)\geq \Zt(H_{b})$ for any $H\in \M$.
\end{cor}

\begin{proof}
Join $H_{b}$ to $H$ by a path \eqref{eqgeod}. By Lemma \ref{lemfirst}, $H_{b}$ is a critical point of the convex function $t\mapsto\Zt(H_{t})$.
\end{proof}

Let $\Aut(X,L)$ be the group of pairs $(f,\tilde{f})$, where $f$ is an automorphism of $X$ and $\tilde{f}:L\To L$ is a holomorphic bundle isomorphism covering $f$. It acts on $W$ by $(\tilde{f}s)(x)=\tilde{f}(s(f^{-1}(x)))$. Since $L$ is very ample, this representation is faithful. We denote its image in $GL(W)$ by $G$. Note that $G$ contains the scalars $\C^{*}$, and that $g\in GL(W)$ belongs to $G$ if and only if the induced automorphism of $\mathbb{P}(W^{*})$ preserves $\iota(X)$, where $\iota:X\To\mathbb{P}(W^{*})$ is the canonical embedding. Indeed, such $g$ induces an automorphism $f$ of $X$, and since $\iota^{*}\mathcal{O}(1)=L$, the natural lift of $g$ to $\mathcal{O}(1)$ restricts to a lift $\tilde{f}$ of $f$ to $L$. In particular, $G$ is a closed subgroup of $GL(W)$. For $g\in G$ corresponding to $(f,\tilde{f})$, $H\in \M$ and $h\in \K$, we define
$$g^{*}H(s,t)=H(gs,gt), \qquad g^{*}h(v,w)=h(\tilde{f}v,\tilde{f}w).$$

\begin{lem}\label{lemequiv}
For any $g\in G$, $H\in \M$ and $h\in \K$, we have
$$\Hilb(g^{*}h)=g^{*}\Hilb(h),\qquad \FS(g^{*}H)=g^{*}\FS(H),$$
$$T(g^{*}H)=g^{*}T(H),\qquad P(g^{*}H)=g^{-1}P(H)g.$$
If in addition $(X,L)$ admits a balanced metric, then $\Zt(g^{*}H)=\Zt(H)$.
\end{lem}

\begin{proof}
For $g \in G$ corresponding to $(f,\tilde{f})$ and $s,t\in W$, we have $(g^{*}h)(s,t)=f^{*}\big(h(gs,gt)\big)$ and $d\mu_{g^{*}h}=f^{*}d\mu_{h}$. Therefore the maps $\Hilb$, $\FS$, $T$ and $P$ are equivariant.

For the last statement, let $H_{t}$ be a smooth path in $\M$ and let $a_{t}$ be the self-adjoint endomorphisms with respect to $H_t$ satisfying $\frac{d}{dt}{H}_{t}(s,t')=H_{t}(a_{t}s,t').$ Then $\frac{d}{dt}(g^*{H}_{t}(s,t'))=g^*H_{t}((g^{-1}a_{t}g)s,t').$ Therfore,  Lemma \ref{lemfirst} gives
$$\frac{d}{dt}\Zt(g^{*}H_{t})=-\frac{V}{N}\tr\big(g^{-1}a_{t}g\,g^{-1}P(H_{t})g\big)=\frac{d}{dt}\Zt(H_{t}).$$
Since $\M$ is connected, $c(g)=\Zt(g^{*}H)-\Zt(H)$ does not depend on $H$. Hence $c:G\To\R$ is a homomorphism. If $H_{b}$ is balanced, Corollary \ref{cormin} gives $c(g)=\Zt(g^{*}H_{b})-\Zt(H_{b})\geq 0$ for all $g\in G$. Applying this to $g^{-1}$, we get $c=0$. Note that the last statement does not hold in
general without the assumption on the existence of balanced metrics(see Appendix A).

\end{proof}

\begin{prop}\label{propB}
Suppose that $H_{b}\in \M$ is balanced. Then the set $\B$ of balanced inner products is the orbit of $H_{b}$ under $G$, i.e.
$$\B=\set{H\in \M\mid T(H)=H}=\set{g^{*}H_{b}\mid g\in G}.$$
Moreover, for any $H\in \B$ there exists $g\in G$, self-adjoint and positive with respect to $H_{b}$, such that $H=g^{*}H_{b}$.
\end{prop}

\begin{proof}
By Lemma \ref{lemequiv}, $g^{*}H_{b}$ is balanced for any $g\in G$. Conversely, let $H$ be balanced. Let $a$ be the $H_{b}$-self-adjoint endomorphism with $H(s,t)=H_{b}(e^{\frac{a}{2}}s,e^{\frac{a}{2}}t)$ and let $H_{t}$ be the path \eqref{eqgeod} from $H_{b}$ to $H$. By Lemma \ref{lemsecond} and Corollary \ref{cormin}, $\Zt(H_{t})$ is convex and attains its minimum at $t=0$ and $t=1$. Hence it is constant on $[0,1]$, and $\frac{d^{2}}{dt^{2}}\big|_{t=0}Z(H_{t})=0$. By Lemma \ref{lemsecond}, $e^{ta}$ preserves $X_{0}=\iota_{H_{b}}(X)$ for all $t$. Therefore $g=e^{\frac{a}{2}}\in G$ and $H=g^{*}H_{b}$.
\end{proof}

As a corollary, we can conclude the analogue of Matsushima's theorem (\cite{Mat}) for balanced metrics.
\begin{cor}\label{correductive}
Assume that $(X,L)$ admits a balanced metric. Then the group $\Aut(X,L)$, which we identify with $G$, is reductive.
\end{cor}

\begin{proof}
     Fix a balanced metric $H_0 \in \M$ and denote by $x^\dagger$ the adjoint of $x\in\textrm{End}(W)$ with respect to $H_0$. Proposition \ref{propB} implies that every element of the orbit $\set{g^*H_0\mid g\in G}$ is of the form $p^*H_0$, where $p\in G$ is self-adjoint and positive with respect to $H_0$. Let $g\in G$, and choose such $p$ with $g^*H_0=p^*H_0$. Then $k=gp^{-1}$ preserves $H_0$, so that $k^{-1}=k^\dagger$. This implies that   $g^\dagger=p^\dagger k^\dagger=pk^{-1}\in G$. Hence $G$ is closed under the adjoint with respect to $H_0$. Now a standard group theoretic argument implies that $G$ is reductive (c.f. \cite{Mo}). Let $R$ be the unipotent radical of $G$. Since the map $x\mapsto(x^\dagger)^{-1}$ is an automorphism of $G$, it preserves $R$ by the uniqueness of $R$.  Therefore, $R$ is also closed under the adjoint. For $u\in R$, the element $u^\dagger u\in R$ is unipotent and positive. Hence $u^\dagger u=I$, and $u$ is unitary and unipotent, i.e. $u=I$. This implies that $R$ is trivial and hence, $G$ is reductive.
\end{proof}
Note that Proposition \ref{propB} is the analogue for balanced metrics of the uniqueness part of the Kempf--Ness theorem. Corollary \ref{correductive} can also be deduced from the reductivity of stabilizers of polystable points, applied to the Chow point of $\iota(X)$ via the work of Zhang and Luo (\cite{Zh}, \cite{L}). 

\begin{rem}\label{rempolar}
The last statement of Proposition \ref{propB} does not hold for every $G$-orbit in $\M$. Indeed by the argument in the proof of Corollary \ref{correductive}, if the last statement of Proposition \ref{propB} holds for the orbit of an element $H \in \M$, then $G$ is closed under the adjoint with respect to $H$ (c.f. Example \ref{excp1}).    
\end{rem}

We now go over the geometry of $\M$. For $H\in \M$, the tangent space $T_{H}\M$ is identified with the space $\Herm_{H}(W)$ of $H$-self-adjoint endomorphisms of $W$ via $a\mapsto H(a\cdot,\cdot)$. We put $\langle a,b\rangle_{H}=\tr(ab)$ and $\norm{a}^2=\tr(a^{2})$. This is the metric of the symmetric space $\M\cong GL(N,\C)/U(N)$. Its exponential map at $H$ is
$$\Exp_H(a)=(e^{a/2})^*H,$$
i.e. $\Exp_H(a)(s,t)=H(e^{a/2}s,e^{a/2}t)$, so that the geodesics through $H$ are the paths \eqref{eqgeod}, and the distance is the one given in the introduction. The space $\M$ is a Hadamard manifold, i.e. it is complete, simply connected and nonpositively curved. Hence $\Exp_H:\Herm_H(W)\to\M$ is a diffeomorphism, any two points of $\M$ are joined by a unique geodesic, and closed bounded subsets of $\M$ are compact. The group $GL(W)$ acts on $\M$ by isometries $H\mapsto g^{*}H$, and the differential of $g^*$ maps $a\in\Herm_H(W)$ to $g^{-1}ag\in\Herm_{g^*H}(W)$. Note that a sequence $H_{m}$ converges to $H$ in $\M$, i.e. the matrices $[H_{m}(s_{i},s_{j})]$ converge for a fixed basis, if and only if $\dist(H_{m},H)\To 0$. By Lemma \ref{lemfirst}, the gradient of $\Zt$ is
\begin{equation}\label{eqgrad}
\nabla\Zt(H)=-\frac{V}{N}P(H), \qquad \norm{\nabla\Zt(H)}_{H}=\frac{V}{N}\norm{P(H)}.
\end{equation}
Moreover, the Hessian of $\Zt$ is given by
$$\nabla^2\Zt(H)(a,a)=\frac{d^2}{dt^2}\Big|_{t=0}\Zt\big(\Exp_H(ta)\big).$$
By Lemma \ref{lemsecond}, it is positive semi-definite for every $H\in\M$.
The proof of the following is standard in the geometry of homogeneous spaces. For completeness, we include the proof.
\begin{prop}\label{propgeomB}
Suppose that $H_b\in\M$ is balanced and let
\begin{align*}
\pp_{H_b}&=\set{a\in\Herm_{H_b}(W)\mid e^{ta}\in G\textrm{ for all }t\in\R}\\
&=\Lie(G)\cap\Herm_{H_b}(W).
\end{align*}
Then $\B=\Exp_{H_b}(\pp_{H_b})$. Moreover, the following hold.
\begin{enumerate}
\item The set $\B$ is a closed, connected, geodesically convex and totally geodesic submanifold of $\M$, and $T_{H_b}\B=\pp_{H_b}$.
\item For every $H\in\M$, there is a unique point $\pi(H)\in\B$ such that $\dist(H,\pi(H))=\dist(H,\B)$. The map $\pi:\M\to\B$ is $G$-equivariant, and $\dist(g^*H,\B)=\dist(H,\B)$ for $g\in G$.
\item Let $\Nor_{H_b}=(T_{H_b}\B)^{\perp}$. For any $y\in\Nor_{H_b}$, we have $\pi(\Exp_{H_b}(y))=H_b$ and $\dist(\Exp_{H_b}(y),\B)=\norm{y}$. Conversely, any $H\in\M$ is of the form $H=\Exp_{\pi(H)}(y)$ with $y\in\Nor_{\pi(H)}$.
\item The quadratic form $\nabla^2\Zt(H_b)$ vanishes exactly on $T_{H_b}\B$. In particular, it is positive definite on $\Nor_{H_b}$.
\end{enumerate}
\end{prop}

\begin{proof}
The set $\pp_{H_b}$ is a real linear subspace of $\Herm_{H_b}(W)$ since $G$ is a closed subgroup og $GL(W)$. For $a\in\pp_{H_b}$, we have  $\Exp_{H_b}(a)=(e^{a/2})^*H_b\in\B$ . Conversely, the proof of Proposition \ref{propB} shows that any $H\in\B$ is of the form $\Exp_{H_b}(a)$, where $e^{ta}$ preserves $\iota_{H_b}(X)$ for all $t$, i.e. $a\in\pp_{H_b}$. Hence $\B=\Exp_{H_b}(\pp_{H_b})$.

Since $\Exp_{H_b}$ is a diffeomorphism , $\B$ is a closed, connected, embedded submanifold and $T_{H_b}\B=\pp_{H_b}$. Let $H_1\in\B$. Hence $H_1=\Exp_{H_b}(a)$ for some $a\in\pp_{H_b}$, and the geodesic $t\mapsto\Exp_{H_b}(ta)$ lies in $\B$. Therefore $\B$ is geodesically convex. Moreover,  the whole geodesic $t\mapsto\Exp_{H_b}(ta)$ lies in $\B$. Hence $\B$ is totally geodesic. This proves $(1)$.

Since $\M$ is a Hadamard manifold, it is a complete $CAT(0)$ space, and $\B$ is a closed convex subset. Hence every $H\in\M$ has a unique nearest point $\pi(H)\in\B$ (c.f.~\cite[Proposition II.2.4]{BH}). By Proposition \ref{propB}, $G$ acts on $\M$ by isometries preserving $\B$. Therefore uniqueness implies $\pi(g^*H)=g^*\pi(H)$ and $\dist(g^*H,\B)=\dist(H,\B)$. This proves $(2)$.

Let $y\in\Nor_{H_b}$, $H=\Exp_{H_b}(y)$ and $H_1\in\B$. By $(1)$, the geodesic from $H_b$ to $H_1$ lies in $\B$, so its initial velocity is tangent to $\B$ and orthogonal to $y$. Since $\M$ is nonpositively curved, the comparison theorem for the geodesic triangle with vertices $H_b$, $H$ and $H_1$ gives
$$\dist(H,H_1)^2\geq\dist(H_b,H)^2+\dist(H_b,H_1)^2=\norm{y}^2+\dist(H_b,H_1)^2.$$
Hence $\dist(H,\B)=\norm{y}$, and the minimum is attained only at $H_1=H_b$. Conversely, let $H\in\M$ and let $y=\Exp_{\pi(H)}^{-1}(H)$. Since $\pi(H)$ minimizes the distance from $H$ to the submanifold $\B$, the first variation formula implies that $y\in\Nor_{\pi(H)}$. This proves $(3)$.

Finally, since $\Zt-Z$ is affine along geodesics,  Lemma \ref{lemsecond} implies $(4)$. 
\end{proof}

\begin{rem}\label{remorbit}

 Proposition \ref{propgeomB} does not hold for an arbitrary $G$-orbit in $\M.$ Indeed, it holds for the orbit $O$ of $H_0 \in \M$ provided that $G$ is closed under the adjoint with respect to $H_0$. Indeed, in this case, we can easily show that $O=\Exp_{H_0}(\pp_{H_0})$ and the same holds at every point of $O$. Since the proof of Proposition \ref{propgeomB} uses the existence of a balanced metric only through the identity $\B=\Exp_{H_1}(\pp_{H_1})$ for $H_1\in\B$, it applies to $O$.
 
 For an arbitrary orbit $O$, the following parts remain true. The orbit $O$ is a closed embedded submanifold of $\M$, since the closed subgroup $G$ of $GL(W)$ acts properly on $\M$, and $\dist(\cdot,O)$ is $G$-invariant. Since $[\pp_{H_0},[\pp_{H_0},\pp_{H_0}]]\subseteq\pp_{H_0}$, the set $\Exp_{H_0}(\pp_{H_0})$ is a closed totally geodesic submanifold contained in $O$. Finally, for every $H_0\in\M$, Lemma \ref{lemsecond} shows as in the proof of Proposition \ref{propgeomB}(4) that $\nabla^2\Zt(H_0)$ is positive semi-definite and vanishes exactly on $\pp_{H_0}$. 
\end{rem}
 In the following example, the orbit of a non-balanced inner product has real dimension $6$, while $\pp_{H_0}$ has real dimension $2$, and the orbit is neither totally geodesic nor geodesically convex. 
\begin{ex}\label{excp1}
Let $X=\CP^1$ and $L=\mathcal{O}(2)$. Then $W$ is the space of quadratic forms $f=a\,z_0^2+b\,z_0z_1+c\,z_1^2$, and $G$ is the image of the representation $\rho$ of $GL(2,\C)$ on $W$ defined by
$$\big(\rho(A)f\big)(z_0,z_1)=f(\alpha z_0+\gamma z_1,\ \beta z_0+\delta z_1),\qquad A=\begin{pmatrix}\alpha&\beta\\ \gamma&\delta\end{pmatrix}.$$
Note that $\rho(AB)=\rho(A)\rho(B)$ and $\rho(\lambda I)=\lambda^2\,\textrm{Id}$, so that $G$ contains the scalars. In the basis $z_0^2,z_0z_1,z_1^2$, the matrix of $\rho(A)$ is
$$M(A)=\begin{pmatrix}\alpha^2&\alpha\beta&\beta^2\\ 2\alpha\gamma&\alpha\delta+\beta\gamma&2\beta\delta\\ \gamma^2&\gamma\delta&\delta^2\end{pmatrix}.$$
The discriminant $\Delta(f)=b^2-4ac$ satisfies $\Delta(\rho(A)f)=\det(A)^2\Delta(f)$, and a nonzero $f$ is a multiple of the square of a linear form if and only if $\Delta(f)=0$. In particular, the elements of $G$ map squares of linear forms to multiples of squares of linear forms. Below, $B^\dagger$ denotes the conjugate transpose of a matrix $B$.

Let $H_b$ be the inner product for which $z_0^2$, $\sqrt2\,z_0z_1$, $z_1^2$ are orthonormal. In this basis, the matrix of $\rho(A)$ is
$$M_b(A)=\begin{pmatrix}\alpha^2&\sqrt2\,\alpha\beta&\beta^2\\ \sqrt2\,\alpha\gamma&\alpha\delta+\beta\gamma&\sqrt2\,\beta\delta\\ \gamma^2&\sqrt2\,\gamma\delta&\delta^2\end{pmatrix},$$
and a direct computation gives $M_b(A)^\dagger=M_b(A^\dagger)$. Hence the adjoint of $\rho(A)$ with respect to $H_b$ is $\rho(A^\dagger)$. In particular, $G$ is closed under the adjoint with respect to $H_b$, and $\rho(A)$ is $H_b$-unitary for $A\in U(2)$. Therefore $H_b$ is invariant under $U(2)$, and it is balanced by the argument of Section \ref{secproj}. Moreover, if $A=UP$ is the polar decomposition of $A$, then $\rho(U)$ is $H_b$-unitary, and $\rho(P)=\rho(P^{1/2})^2$ is self-adjoint and positive with respect to $H_b$. Hence $\rho(A)^*H_b=\rho(P)^*H_b$, in accordance with Proposition \ref{propB}.

Now let $H_0$ be the inner product for which $z_0^2$, $z_0z_1$, $z_1^2$ are orthonormal. Then the adjoint of $\rho(A)$ with respect to $H_0$ has the matrix $M(A)^\dagger$ in this basis. Let $g=\rho(A_1)$, where $A_1=\begin{pmatrix}1&1\\0&1\end{pmatrix}$, i.e. $g$ is the substitution $z_1\mapsto z_0+z_1$. Then
$$M(A_1)=\begin{pmatrix}1&1&1\\0&1&2\\0&0&1\end{pmatrix},\qquad M(A_1)^\dagger=\begin{pmatrix}1&0&0\\1&1&0\\1&2&1\end{pmatrix}.$$
Hence the adjoint $g^\dagger$ of $g$ with respect to $H_0$ is given by
$$g^\dagger(z_0^2)=z_0^2+z_0z_1+z_1^2,\qquad g^\dagger(z_0z_1)=z_0z_1+2z_1^2,\qquad g^\dagger(z_1^2)=z_1^2.$$
Since $\Delta(z_0^2)=0$ and $\Delta(z_0^2+z_0z_1+z_1^2)=-3$, the map $g^\dagger$ does not map squares of linear forms to multiples of squares of linear forms. Therefore $g^\dagger\notin G$, i.e. $G$ is not closed under the adjoint with respect to $H_0$. In particular, $H_0$ is not balanced, by Proposition \ref{propB} and Remark \ref{rempolar}. Note that with respect to $H_b$, the adjoint of $g$ is $\rho(A_1^\dagger)$, i.e. the substitution $z_0\mapsto z_0+z_1$, which lies in $G$.

The element $g^*H_0$ in the basis $z_0^2,z_0z_1,z_1^2$ is 
$$M(A_1)^\dagger M(A_1)=\begin{pmatrix}1&1&1\\1&2&3\\1&3&6\end{pmatrix}.$$
 Suppose that $g^*H_0=p^*H_0$ for some $p\in G$ which is self-adjoint and positive with respect to $H_0$. Then $k=gp^{-1}$ preserves $H_0$, and therefore $g^\dagger=pk^{-1}\in G$, which is a contradiction. Hence $g^*H_0$ is not of the form $p^*H_0$, and the last statement of Proposition \ref{propB} does not hold for the orbit of $H_0$.

Finally, we describe the orbit $O$ of $H_0$ near $H_0$. Using the basis $z_0^2,z_0z_1,z_1^2$, we identify $\Herm_{H_0}(W)$ with the space of Hermitian $3\times3$ matrices. Differentiating $M(A)$ at $A=I$, we see that $\Lie(G)$ is spanned over $\C$ by
\begin{align*}
D_0&=\begin{pmatrix}2&0&0\\0&1&0\\0&0&0\end{pmatrix},& D_1&=\begin{pmatrix}0&0&0\\0&1&0\\0&0&2\end{pmatrix},\\
E&=\begin{pmatrix}0&0&0\\2&0&0\\0&1&0\end{pmatrix},& F&=\begin{pmatrix}0&1&0\\0&0&2\\0&0&0\end{pmatrix}.
\end{align*}
The velocity of the path $t\mapsto(e^{tx})^*H_0$ at $t=0$ is $x+x^\dagger$. Writing $x=aD_0+bD_1+cE+dF$ with $a,b,c,d\in\C$, we get
$$x+x^\dagger=2\,\textrm{Re}(a)D_0+2\,\textrm{Re}(b)D_1+\begin{pmatrix}0&2\bar c+d&0\\2c+\bar d&0&\bar c+2d\\0&c+2\bar d&0\end{pmatrix}.$$
Hence
$$T_{H_0}O=\set{Y\mid Y^\dagger=Y,\ Y_{13}=0,\ 2Y_{22}=Y_{11}+Y_{33}},$$
which has real dimension $6$. On the other hand, $x$ is Hermitian if and only if $a,b\in\R$, $2c=\bar d$ and $c=2\bar d$, i.e. $a,b\in\R$ and $c=d=0$. Hence $\pp_{H_0}=\set{rD_0+sD_1\mid r,s\in\R}$ has real dimension $2$, and $\Exp_{H_0}(\pp_{H_0})$ is a proper subset of $O$. Moreover, let $E_{ij}$ denote the elementary matrices. Then
$$u=E_{12}+E_{21},\qquad v=E_{23}+E_{32},\qquad w=D_0$$
lie in $T_{H_0}O$, while
$$\big[[u,v],w\big]=\big[E_{13}-E_{31},D_0\big]=-2(E_{13}+E_{31})\notin T_{H_0}O.$$
Hence $O$ is not totally geodesic. Since a closed geodesically convex submanifold is totally geodesic, $O$ is not geodesically convex either. For comparison, the orbit $\B$ of $H_b$ has real dimension $4$, and $T_{H_b}\B=\pp_{H_b}$ by Proposition \ref{propgeomB}. Note that this agrees with dimension count. Indeed, $G$ has real dimension $8$. The stabilizer of $H_b$ in $G$ is the image of $U(2)$ which has real dimension $4$ and the stabilizer of $H_0$ is the image of diagonal unitary matices whoes dim is $2.$
\end{ex}

Note that if $\Aut(X,L)/\C^{*}$ is discrete, then $\pp_{H_b}=\R\cdot I$ and $\B=\set{e^{\alpha}H_b\mid\alpha\in\R}$. In this case, Proposition \ref{propB} is Donaldson's uniqueness theorem (\cite{D1}).

Throughout the rest of the paper, we assume that $(X,L)$ admits a balanced metric. We denote the set of balanced inner products by $\B$, and we put $m_{*}=\min_{\M}\Zt$, which is the value of $\Zt$ on $\B$ by Corollary \ref{cormin}. To avoid trivialities, we also assume that $\B\neq\M$. Note that if $\B=\M$, then every $H\in\M$ is balanced, i.e. $T$ is the identity map of $\M$. This happens for $(\CP^n,\mathcal{O}(1))$ (see Section \ref{secproj}).

\section{Convergence to the set of balanced metrics}
The main goal of this section is to prove that the sequence $H_m=T^m(H_0)$ approaches the submanifold $\B$ (Proposition \ref{propconv}). If $\Aut(X,L)/\C^*$ is discrete, Sano uses the properness of $\Zt$ modulo scaling (\cite{S}). In general, $\Zt$ is constant on $\B$, which is noncompact even modulo scaling. Therefore $\Zt$ is not proper, even modulo scaling. However, the existence of balanced metrics implies that the Chow point is poly-stable, and therefore $\Zt$ is ``relatively proper'' by the Kempf--Ness theorem. Here, we only use the convexity property to prove the following properness of $\Zt$ modulo $G$.

\begin{lem}\label{lemproper}
There exist a continuous nondecreasing function $\psi:[0,\infty)\to[0,\infty)$ and a constant $\eta>0$ such that $\psi(0)=0$, $\psi(t)>0$ for $t>0$, $\psi(t)\geq\eta t$ for $t\geq 1$, and
$$\Zt(H)-m_*\geq\psi\big(\dist(H,\B)\big)\qquad\textrm{for all } H\in\M.$$
\end{lem}

\begin{proof}
Fix $H_b\in\B$ and let $S$ be the unit sphere of $\Nor_{H_b}$. For $y\in S$, put
$$f_y(t)=\Zt\big(\Exp_{H_b}(ty)\big)-m_*.$$
The function $f_y$ is convex, $f_y(0)=0$ and $f_y'(0)=0$ since $H_b$ is balanced (Lemma \ref{lemfirst}). Moreover, $f_y''(0)=\nabla^2\Zt(H_b)(y,y)>0$ by Proposition \ref{propgeomB}. Hence convexity implies that $f_y'>0$ on $(0,\infty)$. Therefore $f_y$ is positive and increasing on $(0,\infty)$. By convexity, we also have $f_y(t)\geq tf_y(1)$ for $t\geq1$. Put
$$\psi(t)=\min_{y\in S}f_y(t),\qquad \eta=\psi(1).$$
Since $(t,y)\mapsto f_y(t)$ is continuous and $S$ is compact, $\psi$ has the stated properties. Now let $H\in\M$ such that $\pi(H)=H_b.$ Proposition \ref{propgeomB} implies that there exists $y \in \Nor_{H_b}$ such that $H=\Exp_{H_b}(t_Hy)$, where $t_H=\dist(H, \B).$ Hence,
$$\Zt\big(\Exp_{H_b}(t_Hy)\big)-m_*=f_y(t_H) \geq \psi(t_H)=\psi(\dist(H,\B)).$$

For an arbitrary $H \in \M$, by Proposition \ref{propB}, there exists $g\in G$ with $g^*\pi(H)=H_b$. Since $\Zt$ and $\dist(\cdot,\B)$ are $G$-invariant and $\pi$ is $G$-equivariant, we conclude the lemma. 
\end{proof}

\begin{rem}\label{remKN}
Lemma \ref{lemproper} can also be deduced from the Kempf--Ness theorem (\cite{KN}, c.f. \cite{GRS}). Indeed, fix $H_b\in\B$. By the work of Zhang (\cite{Zh}, see also \cite{L}, \cite{PS}), for $g\in SL(W)$ the difference $\Zt(g^*H_b)-\Zt(H_b)$ is a positive multiple, depending only on $(X,L)$, of the logarithm of the ratio of the norms of $g\cdot\textrm{Chow}(\iota_{H_b}(X))$ and $\textrm{Chow}(\iota_{H_b}(X))$, where the norm is the one induced by $H_b$. Since a balanced metric exists, the Chow point is polystable, and the Kempf--Ness theorem implies that this function is proper modulo the stabilizer of the Chow point in $SL(W)$, which is contained in $G$. Since $\Zt$ is invariant under scaling, we conclude that the sublevel sets of $\Zt$ are compact modulo $G$. The advantage of Lemma \ref{lemproper} is that it only uses the convexity of $\Zt$, and it also gives the function $\psi$ explicitly.
\end{rem}

\begin{prop}\label{propconv}
Let $H_{0}\in \M$ and $H_{m}=T^{m}(H_{0})$. Then
\begin{enumerate}
\item $\log\det H_{m}$ converges to a real number.
\item $P(H_{m})\To 0$.
\item $\Zt(H_{m})$ decreases to $m_*$.
\item $\dist(H_{m},\B)\To 0$.
\end{enumerate}
\end{prop}

\begin{proof}
By Lemma \ref{lemmono}, the sequences $Z(H_{m})$, $\log\det H_{m}$ and $\Zt(H_m)$ are nonincreasing. By Corollary \ref{cormin}, $\Zt\geq m_*$. 
Therefore
$$\frac{V}{N}\log\det H_m=\Zt(H_m)-Z(H_m)\geq m_*-Z(H_0),$$
i.e. $\log\det H_m$ is bounded below. Hence it converges to a real number. This proves $(1)$. By Lemma \ref{lemdescent},
$$\norm{P(H_m)}^2\leq 2N\big(\log\det H_m-\log\det H_{m+1}\big),$$
which tends to zero by $(1)$. This proves $(2)$.

Fix $H_b\in\B$ and let $\psi$ and $\eta$ be as in Lemma \ref{lemproper}. Since $\Zt(H_m)\leq\Zt(H_0)$, Lemma \ref{lemproper} gives $\psi\big(\dist(H_m,\B)\big)\leq\Zt(H_0)-m_*$. Since $\psi(t)\geq\eta t$ for $t\geq1$, we get $\dist(H_m,\B)\leq R_1=\max\big(1,(\Zt(H_0)-m_*)/\eta\big)$ for all $m$. Let $g_m\in G$ be such that $g_m^*\pi(H_m)=H_b$ and let $K_m=g_m^*H_m$. Then $\pi(K_m)=H_b$ and $$\dist(K_m,H_b)=\dist(H_m,\B)\leq R_1.$$ Hence the sequence $\{K_m\}$ is bounded. Moreover, by Lemma \ref{lemequiv} and part $(2)$, we have $$\norm{P(K_m)}=\norm{g_m^{-1}P(H_m)g_m}=\norm{P(H_m)}\to0.$$ By continuity of $P$, every limit point $K_\infty$ of the sequence $K_m$ satisfies $P(K_\infty)=0$, which means $K_\infty\in\B$. Therefore, we conclude that $\dist(K_m,\B)\to0$ and $\Zt(K_m)\to m_*$. Finally, $\dist(K_m,\B)=\dist(H_m,\B)$ and $\Zt(K_m)=\Zt(H_m)$ by Proposition \ref{propgeomB} and Lemma \ref{lemequiv}. This proves $(3)$ and $(4)$.
\end{proof}

Note that the proof of Proposition \ref{propconv} does not use the boundedness of the sequence $H_{m}$, which is the main step in \cite{S}. If $\Aut(X,L)/\C^{*}$ is discrete, then $\B=\set{e^{\alpha}H_{b}\mid\alpha\in\R}$, and $(1)$ and $(4)$ already imply the convergence of $H_{m}$.
Indeed, let $\pi(H_m)=e^{\alpha_m}H_b.$
We have
$$\abs{\log\det H_m-N\alpha_m-\log\det H_b} \leq \sqrt{N}\dist(H_m,\B) \to 0.$$
Hence $\alpha_m$ converges by $(1)$, and so does $H_m$.
\section{Proof of Theorem \ref{mainthm}}

The goal of this section is to prove the main theorem. We give two proofs. The first proof uses the convexity of the functional $\Zt$. Roughly speaking, we use the fact that $\Zt$ is strictly convex in the directions normal to $\B$. It proves the convergence of $T^m(H_0)$ with an exponential rate that does not depend on $H_0$. However, in the first proof, the constant $C=C(H_0)$ in Theorem \ref{mainthm} depends on $H_0$. The second proof relies on the fact that the linearization of $T$ at a balanced metric is self-adjoint and has eigenvalues in $(0,1]$, where the eigenvalue $1$ only arises in the directions tangent to $\B$. It has the advantage of showing the constant $C$ depends only on $\dist(H_0,\B)$. Another advantage of the second proof is that it gives the optimal rate of convergence. Therefore, the second proof gives the full statement of Theorem \ref{mainthm}. At the end of the section, we show that the rate $\rho_*$ is optimal (Proposition \ref{propsharp}).

\subsection{First proof}

We consider a toy model before going over the actual proof. This example illustrates the basic idea behind the first proof. Let $F:\R^2 \to \R$ be a convex $C^2$ function such that $0=F(x,0)\leq F$ and $F_{yy}(x,0)\geq 2c>0$. Note that $F_{y}(x,0)=0$, since $y=0$ is a minimum point of $F(x,\cdot)$. The $C^2$ assumption implies that for any $x$, there exists $\delta_x>0$ such that $F_{yy}(x,y) \geq c$ for all $|y| \leq \delta_x$. For $0\leq y \leq \delta_x$, we have
$$F_{y}(x,y)=F_{y}(x,y)-F_{y}(x,0)=\int_{0}^yF_{yy}(x,t)\,dt\geq cy.$$
Since $F_y(x,\cdot)$ is nondecreasing by convexity, we get
$$F(x,y)=\int_{0}^y F_y(x,t)\,dt \leq yF_{y}(x,y) \leq c^{-1}F_y(x,y)^2.$$
Motivated by this toy example, we prove the following lemma which is crucial in our first proof. We also note that $G$ acts on $\M$ by isometries and it preserves $\Zt$. This simple observation enables us to show that all estimates are uniform along $\B$.

\begin{lem}\label{lemloj}
There exist $\delta>0$ and $c_{1}>0$ such that for any $H\in \M$ with $\dist(H,\B)<\delta$,
$$\norm{\nabla\Zt(H)}_{H}\geq c_{1}\,\dist(H,\B) \qquad \textrm{and} \qquad \Zt(H)-m_*\leq \frac{1}{c_{1}}\norm{\nabla\Zt(H)}_{H}^{2}.$$
\end{lem}

\begin{proof}
Fix $H_b\in\B$ and let $S$ be the unit sphere of $\Nor_{H_b}$. By Proposition \ref{propgeomB}, there exists $c>0$ such that $\nabla^2\Zt(H_b)(y,y)\geq 2c\norm{y}^2$ for any $y\in\Nor_{H_b}$. The map
$$(t,y)\longmapsto\frac{d^2}{dt^2}\Zt\big(\Exp_{H_b}(ty)\big)$$
is continuous on $[0,1]\times S$ and is at least $2c$ at $t=0$. By uniform continuity, there exists $\delta\in(0,1]$ such that it is at least $c$ on $[0,\delta]\times S$.

Let $H\in\M$ with $t_H:=\dist(H,\B)<\delta$.
Choose $g\in G$ with $g^*\pi(H)=H_b$. By Lemma \ref{lemequiv}, we have $\Zt(g^*H)=\Zt(H)$ and $\norm{\nabla\Zt(g^*H)}=\norm{\nabla\Zt(H)}$, since $P(g^*H)=g^{-1}P(H)g$. Moreover, $\dist(g^*H,\B)=t_H$ by Proposition \ref{propgeomB}.  Hence without loss of generality, we may assume that $\pi(H)=H_b$. Then $H=\Exp_{H_b}(t_Hy)$ for some $y\in S$, by Proposition \ref{propgeomB}, where $S$ is the unit sphere in $\Nor_{H_b}$. Define $$\phi(t)=\Zt(\Exp_{H_b}(ty))\,\,\,\, \textrm{for}\,\,\,\, t\in[0,t_H]$$Then $\phi'(0)=0$, since $H_b$ is a critical point, and $\phi''\geq c$ on $[0,t_H]$. Since the geodesic $t\mapsto\Exp_{H_b}(ty)$ has unit speed, we have
$$\norm{\nabla\Zt(H)}_{H} \geq \phi^{\prime}(t_H)=\int_{0}^{t_H}\phi^{\prime\prime}(s)\,ds \geq c\, t_H=c\,\dist(H, \B).$$
Moreover, $\phi'$ is nondecreasing, and therefore
$$\Zt(H)-m_*=\int_{0}^{t_H}\phi^{\prime}(s)\,ds\leq t_H \phi^{\prime}(t_H)\leq \frac{1}{c}\phi'(t_H)^2\leq\frac{1}{c}\norm{\nabla\Zt(H)}^2_{H}.$$
\end{proof}

\begin{proof}[First proof of the convergence in Theorem \ref{mainthm}]
Let $H_m=T^m(H_0)$ and $e_m=\Zt(H_m)-m_*$.
The proof will be divided into several steps:
	
\medskip
\noindent {\bf STEP I (Exponential decay of $e_m$):} 
\noindent
By Proposition \ref{propconv}, there exists $m_0$ such that for $m\geq m_0$, we have $\dist(H_m,\B)<\delta$ and all eigenvalues of $P(H_m)$ have absolute value at most $\frac12$. By \eqref{eqgrad}, Lemma \ref{lemdescent} and Lemma \ref{lemloj}, we have, for $m\geq m_0$,
$$e_m-e_{m+1}\geq\frac{V}{2N^2}\norm{P(H_m)}^2=\frac{1}{2V}\norm{\nabla\Zt(H_m)}_{H_m}^2\geq\frac{c_1}{2V}\,e_m.$$
Therefore $e_{m+1}\leq(1-\kappa)e_m$ for $m\geq m_0$, for some positive $\kappa=\kappa(c_1).$  Since $e_m\leq e_0$ for all $m$, we get $$e_m\leq C(1-\kappa)^{m}\qquad\textrm{for all } m\geq0,$$ where $C=e_0(1-\kappa)^{-m_0}$.

\medskip
\noindent {\bf STEP II (Exponential decay of $\dist(H_m,H_{m+1})$):} 
\noindent
 For $m\geq m_0$, Lemma \ref{lemdescent} implies that
$$
\dist(H_m,H_{m+1})\leq\sqrt2\norm{P(H_m)}\leq\frac{2N}{\sqrt V}\,(e_m-e_{m+1})^{\frac12}
\leq C(1-\kappa)^{\frac{m}{2}}.
$$
Hence $H_m$ is a Cauchy sequence. Let $H_\infty$ be its limit. By continuity of $P$ and Proposition \ref{propconv}, $P(H_\infty)=0$ which implies that $H_\infty$ is balanced. Moreover,
$$\dist(H_m,H_\infty)\leq\sum_{j\geq m}\dist(H_j,H_{j+1})\leq C\lambda_1^m,\qquad \lambda_1=(1-\kappa)^{\frac12}.$$
Therefore, this holds for all $m$ after choosing $C$ larger.  Note that $\lambda_1$ depends only on $(X,L)$, while $C$ depends on $H_0$.
\end{proof}

\subsection{Second proof}\label{subsecsecond}

We start with the derivative of $T$ at a balanced metric. In the toy example of \cite[Section 2.1.1]{D3}, Donaldson observes that, for almost all initial conditions, the asymptotic rate of convergence of the iteration is the largest eigenvalue of the derivative of $T$ at the fixed point. In our setting, $d_HT$ is the identity on $T_H\B$, so the relevant quantity is the largest eigenvalue on $\Nor_H$ (Theorem \ref{thmlin}(5) below). In \cite[Section 4.1]{D3}, Donaldson considers the analogous map $T_\nu$, in which a measure $\nu$ on $X$ is fixed. He shows that its derivative at a fixed point is given by an operator on Hermitian matrices, denoted $Q$ in \cite{D3}. If $H$ is balanced and $\nu=d\mu_{\FS(H)}$, then $H$ is also a fixed point of $T_\nu$, and in the notation of Remark \ref{remberezin} below its derivative is $Q\circ\sigma$. Note that our $Q$, which maps functions to matrices, is not the operator $Q$ of \cite{D3}. In \cite[Corollary 4.9]{IP}, Ioos and Polterovich express the derivative of $T$ itself in terms of the Berezin--Toeplitz quantization and its adjoint. For high powers of a polarization admitting a constant scalar curvature K\"ahler metric with $\Aut(X,L)/\C^*$ discrete, they compute the asymptotics of its spectral gap (\cite[Theorem 1.2(2)]{IP}); see also \cite{LM} and Section \ref{secasym}. For the convenience of the reader, we give a self-contained proof of the following description of $d_HT$, based on Section $2$.

For $H\in\M$, $h=\FS(H)$ and $a\in\Herm_H(W)$, we define the function $\hat a$ on $X$ by
$$\hat a(x)=\sum_{i,j=1}^N H(as_j,s_i)\,h\big(s_i(x),s_j(x)\big),$$
where $s_1,\dots,s_N$ is an $H$-orthonormal basis of $W$. Note that $\hat a$ does not depend on the choice of the basis, and that $\hat a=\mu_a\circ\iota_H$, where $\mu_a$ is the function of Lemma \ref{lemsecond}.

\begin{thm}\label{thmlin}
Let $H\in\B$ and $\Nor_H=(T_H\B)^\perp$. Then the following hold.
\begin{enumerate}
\item $d_HT=I-\frac{N}{V}\nabla^2\Zt(H)$. In particular, $d_HT$ is self-adjoint with respect to $\langle\cdot,\cdot\rangle_H$.
\item $d_HT$ is the identity on $T_H\B$ and maps $\Nor_H$ to itself.
\item For any $a\in\Herm_H(W)$, we have
\begin{equation}\label{eqdT}
\langle d_HT(a),a\rangle_H=\frac{N}{V}\int_X\big(\hat a^2+|\partial\hat a|^2_{\omega_h}\big)\,d\mu_h,
\end{equation}
where $h=\FS(H)$ and $|\partial f|^2_{\omega_h}d\mu_h=\sqrt{-1}\,\partial f\wedge\bar\partial f\wedge\omega_h^{n-1}/(n-1)!$ for real functions $f$.
\item All eigenvalues of $d_HT$ lie in $(0,1]$. Moreover, all eigenvalues of the restriction of $d_HT$ to $\Nor_H$ lie in $(0,1)$.
\item The largest eigenvalue $\rho_*$ of the restriction of $d_HT$ to $\Nor_H$ does not depend on $H\in\B$. Moreover, $\norm{d_HT(y)}\leq\rho_*\norm{y}$ for any $y\in\Nor_H$.
\end{enumerate}
\end{thm}

\begin{proof}
Note that $\M$ can be regarded as an open subset of the real vector space of Hermitian forms on $W$. Let $H_t$ be a smooth path with $H_t=H$ at $t=0$ and initial velocity $a$. Then $T(H_t)=H_t\big((I+P(H_t))\cdot,\cdot\big)$. Differentiating at $t=0$, we get
$$d_HT(a)=a+d_HP(a).$$  
Here we identify tangent vectors at $H$ with $H$-self-adjoint endomorphisms. On the other hand, the vector field $\nabla\Zt=-\frac{V}{N}P$ vanishes at $H$ since $H$ is balanced. Thus its covariant derivative at $H$ agrees with its ordinary derivative in this chart, i.e. $\nabla^2\Zt(H)(a)=-\frac{V}{N}d_HP(a)$. Therefore $d_HT=I-\frac{N}{V}\nabla^2\Zt(H)$. This proves $(1)$ and $(2)$.
For $(3)$, let $s_1,\dots,s_N$ be an $H$-orthonormal basis with $as_i=\lambda_is_i$, and let $H_t$ be the path \eqref{eqgeod}. By \eqref{eqphi}, $\FS(H_t)=e^{-\varphi_t}h$ with $\dot\varphi_0=-\sum\lambda_i|s_i|_h^2=-\hat a$. Since $\omega_{e^{-\varphi_t}h}=\omega_h+\ddbar\varphi_t$, we have, for $s,s'\in W$,
$$\frac{d}{dt}\Big|_{t=0}T(H_t)(s,s')=\frac{N}{V}\int_X h(s,s')\big(\hat a-\Delta_h\hat a\big)\,d\mu_h,$$
where $\Delta_hf\,d\mu_h=\ddbar f\wedge\omega_h^{n-1}/(n-1)!$. The left hand side is $H(d_HT(a)s,s')$. Hence, for any $b\in\Herm_H(W)$,
\begin{align*}
\tr\big(d_HT(a)\,b\big)&=\sum_{i,j=1}^N H(bs_j,s_i)\,H\big(d_HT(a)s_i,s_j\big)\\
&=\frac{N}{V}\int_X\big(\hat a-\Delta_h\hat a\big)\hat b\,d\mu_h.
\end{align*}
Taking $b=a$ and integrating by parts, we get \eqref{eqdT}. Note that this agrees with $(1)$ and the proof of Lemma \ref{lemsecond}. Indeed, since $H$ is balanced, we have $\int_X\mu_{a^2}\circ\iota_H\,d\mu_h=\frac{V}{N}\sum\lambda_i^2T(H)(s_i,s_i)=\frac{V}{N}\norm{a}^2$, and therefore the proof of Lemma \ref{lemsecond} gives
\begin{align*}
\nabla^2\Zt(H)(a,a)&=\frac{V}{N}\norm{a}^2-\int_X\big(\hat a^2+|\partial\hat a|^2_{\omega_h}\big)\,d\mu_h\\
&=\frac{V}{N}\Big(\norm{a}^2-\langle d_HT(a),a\rangle_H\Big).
\end{align*}

Now, $(4)$ and $(5)$ follow from 
\eqref{eqdT}. Note that $G$ acts transitively on $\B$ by isometries and $T$ is $G$-equivariant. Therefore, the spectrum of $d_HT$ does not depend on $H \in \B.$
\end{proof}

\begin{rem}\label{remberezin}
The proof of Theorem \ref{thmlin}(3) gives a useful description of $d_HT$. Let $H\in\B$, $h=\FS(H)$, and let $s_1,\dots,s_N$ be an $H$-orthonormal basis. We have \[f \in C^{\infty}(X,\R) \to Q(f) \in \Herm_H(W)\]  defined by
$$H\big(Q(f)s,s'\big)=\frac{N}{V}\int_X f\,h(s,s')\,d\mu_h.$$
Put $\sigma(a)=\hat a$. Then $\tr(Q(f)b)=\frac{N}{V}\int_Xf\,\sigma(b)\,d\mu_h$, i.e. $\sigma$ is the adjoint of $Q$. The proof of Theorem \ref{thmlin}(3) shows that
$$d_HT=Q\circ(1-\Delta_h)\circ\sigma.$$
After some normalization, this is the formula $D_q\mathcal{T}_p=T_p\big(1+\frac{1}{4\pi}\Delta_{\FS(q)}\big)T_p^*$ of \cite[Corollary 4.9]{IP},  where $\Delta_{\FS(q)}$ is the positive Laplacian.
The composition $\mathfrak{B}=\sigma\circ Q$ is the Berezin transform of $h$, given by
$$\mathfrak{B}f(x)=\frac{N}{V}\int_X\Big(\sum_{i,j=1}^N h\big(s_i(x),s_j(x)\big)\,h\big(s_j(y),s_i(y)\big)\Big)f(y)\,d\mu_h(y).$$
The integral kernel is the squared pointwise norm of the Bergman kernel of $H$. At a balanced metric, $\mathfrak{B}$ coincides with the operator denoted $Q_K$ in \cite{D3} and \cite{LM}. Write $d_HT=YY^*$ with $Y=Q(1-\Delta_h)^{1/2}$. Since $d_HT$ is injective by Theorem \ref{thmlin}(4), its eigenvalues, counted with multiplicity, are exactly the nonzero eigenvalues of the self-adjoint operator $Y^*Y=(1-\Delta_h)^{1/2}\,\mathfrak{B}\,(1-\Delta_h)^{1/2}$ on $L^2(X,\frac{N}{V}d\mu_h)$.
\end{rem}

\begin{lem}\label{lemlocal}
Fix $\lambda$ with $\rho_{*}<\lambda<1$. There exist $r_0\in(0,1]$ and $C_{1}>0$ such that for any $H\in\M$ with $\dist(H,\B)<r_0$, we have
\begin{enumerate}
\item $\dist\big(T(H),\B\big)\leq \rho_*\dist(H,\B)+C_1\dist(H,\B)^2\leq \lambda\,\dist(H,\B)$,
\item $\dist\big(H,T(H)\big)\leq C_{1}\,\dist(H,\B)$.
\end{enumerate}
\end{lem}

\begin{proof}
Fix $H_b\in\B$. As in the proof of Lemma \ref{lemloj}, both sides of $(1)$ and $(2)$ are $G$-invariant, by Lemma \ref{lemequiv} and Proposition \ref{propgeomB}. Hence we may assume that $\pi(H)=H_b$, i.e. $H=\Exp_{H_b}(y)$ with $y\in\Nor_{H_b}$ and $\norm{y}=\dist(H,\B)$. Define
$$\Phi=\Exp_{H_b}^{-1}\circ T\circ\Exp_{H_b}:\Herm_{H_b}(W)\To\Herm_{H_b}(W).$$

Note that $\Phi$ is well-defined since the exponential map is a global diffeomorphism from $T_{H_b}\M$ onto $\M$.
 We have $\Phi(0)=0$ and $d_0\Phi=d_{H_b}T$. Since the map $\Phi$ is smooth, there exists $A_1>0$ such that
$$\norm{\Phi(y)-d_{H_b}T(y)}\leq A_1\norm{y}^2\qquad\textrm{for }\norm{y}\leq1.$$
Applying Theorem \ref{thmlin}, we have
\[\norm{\Phi(y)}\leq(1+A_1)\norm{y}, \qquad y \in \Nor_H\,\,\textrm{with }\norm{y}\leq1 .\]  Since $H_b\in\B$, we have
\begin{align*}
\dist(T(H),\B)&\leq\dist\big(\Exp_{H_b}(\Phi(y)),\Exp_{H_b}(d_{H_b}T(y))\big)+\norm{d_{H_b}T(y)}\\
&\leq A_2\norm{\Phi(y)-d_{H_b}T(y)}+\rho_*\norm{y}\\&\leq A_1A_2\norm{y}^2+\rho_*\norm{y},
\end{align*}
where $A_2=\norm{\nabla \Exp_{H_b}}_{B_{(1+A_1)}(0)}.$
This implies $(1)$.
To prove $(2),$ we have
\begin{align*}
\dist(H,T(H))&=\dist\big(\Exp_{H_b}(y),\Exp_{H_b}(\Phi(y))\big)\\
&\leq A_2\norm{y-\Phi(y)}\leq A_2(\norm{y}+\norm{\Phi(y)})\\&\leq C_1\norm{y}.
\end{align*}

\end{proof}

\begin{lem}\label{lementer}
Let $r_0$ be as in Lemma \ref{lemlocal}. For any $R>0$, there exists an integer $m_0=m_0(R)$ such that $\dist(T^m(H),\B)<r_0$ for all $m\geq m_0$ and all $H\in\M$ with $\dist(H,\B)\leq R$.
\end{lem}

\begin{proof}
Fix $H_b\in\B$ and let $K_R$ be the closed ball of radius $R$ about $H_b$, which is compact. For any $H\in K_R$, Proposition \ref{propconv} gives an integer $m(H)$ with $\dist(T^{m(H)}(H),\B)<r_0$. By continuity, the same holds on an open neighbourhood $U_H$ of $H$. By compactness, $K_R$ is covered by finitely many $U_{H_1},\dots,U_{H_k}$. Now we choose $m_0$ as the maximum of $m(H_i)$. Applying Lemma \ref{lemlocal} concludes the statement.  Note that for $H \in \M$ with $\dist(H,\B)\leq R$, there exists $g\in G$ with $g^*\pi(H)=H_b$, so that $g^*H\in K_R$. 
\end{proof}

\begin{proof}[Proof of Theorem \ref{mainthm}]
If $T$ is the identity map, there is nothing to prove. Hence we may assume that $\B\neq\M$. Fix $R>0$ and  $\lambda\in(\rho_*,1)$. Let $r_0$ and $C_1$ be as in Lemma \ref{lemlocal} and $m_0=m_0(R)$ be as in Lemma \ref{lementer}. Let $H_0\in\M$ with $\dist(H_0,\B)\leq R$, and put $H_m=T^m(H_0)$ and $d_m=\dist(H_m,\B)$.

\medskip
\noindent {\bf STEP I (Exponential convergence of $H_m$):}
\noindent
 Lemma \ref{lementer} implies that $d_m<r_0$ for $m\geq m_0$. Hence for $m \geq m_0$, we have $d_{m+1}\leq\lambda d_m$ by Lemma \ref{lemlocal}, and therefore
$$d_m\leq d_{m_0}\lambda^{m-m_0}\leq r_0\lambda^{m-m_0} \leq \lambda^{m-m_0}.$$
By Lemma \ref{lemlocal} again,
$$\dist(H_m,H_{m+1})\leq C_1d_m\leq C_1\lambda^{m-m_0}.$$
Hence $H_m$ is a Cauchy sequence and therefore, it converges to a limit $H_{\infty}.$ Since $d_m\to0$ and $\B$ is closed, $H_\infty\in\B$.

\medskip
\noindent {\bf STEP II (Convergence rate improvement):}
\noindent
We now improve the rate from $\lambda$ to $\rho_*$. By Lemma \ref{lemlocal}(1), we have $d_{m+1}\leq(\rho_*+C_1d_m)d_m$ for $m\geq m_0$. 
Therfore, \[\log d_m \leq \log \rho_*+\log (1+C_1\rho_*^{-1}d_{m-1})+\log d_{m-1}.\]
Therefore, for $m\geq m_0$, we get
\begin{align*}\log d_m &\leq \log d_{m_0} + (m-m_0)\log \rho_*+C_1\rho_*^{-1}\sum_{j=m_0}^{m-1}d_j\\& \leq (m-m_0)\log \rho_*+C_1\rho_*^{-1}\sum_{j=m_0}^{\infty}d_j\\& \leq (m-m_0)\log \rho_*+C_2.
\end{align*}
 Therefore, for any $m \geq m_0$, we have
$$\dist(H_m,H_\infty)\leq\sum_{j=m}^{\infty}\dist(H_j,H_{j+1})\leq C_1\sum_{j=m}^{\infty}d_j\leq C_3\rho_*^m.$$

\medskip
\noindent {\bf STEP III ($C=C(R)$):}
\noindent
We have already proved that there exist $m_0$ and $C_2$, depending only on $R$ and $(X,L)$, such that
$$\dist(H_m,H_{\infty}) \leq C_2\rho_*^m\qquad\textrm{for } m\geq m_0.$$
It remains to show that $\dist(H_m,H_\infty)$ is uniformly bounded for $m<m_0$ in terms of $R$. Let $K_R$ be as in the proof of Lemma \ref{lementer} and $A_R=\max_{K_R}\Zt-m_*$. Note that by the $G$-invariance of $\Zt$, $A_R$ only depends on $R$ and $(X,L)$. As before, we may assume that $H_0\in K_R$. Hence $\Zt(H_j)\leq\Zt(H_0)\leq m_*+A_R$ for all $j$ by Lemma \ref{lemmono}. Therefore Lemma \ref{lemproper} implies that $\dist(H_j,\B)\leq R_1$ for all $j$, where $R_1=\max(1,A_R/\eta)$. The function $H\mapsto\dist(H,T(H))$ is continuous and $G$-invariant. Hence it is bounded by $D=\max_{K_{R_1}}\dist(\cdot,T(\cdot))$ on the set $\set{H\mid\dist(H,\B)\leq R_1}$. Therefore, for $m<m_0$,
$$\dist(H_m,H_\infty)\leq\sum_{j=m}^{m_0-1}\dist(H_j,H_{j+1})+\dist(H_{m_0},H_\infty)\leq m_0D+C_3 \rho_*^{m_0},$$
which only depends on $R$ and $(X,L).$ 
\end{proof}

\begin{rem}\label{remC}
The constant $C$ cannot be chosen independent of $R$. Indeed, taking $m=0$ in Theorem \ref{mainthm}, we get $C\geq\dist(H_0,H_\infty)\geq\dist(H_0,\B)$, and $\dist(H_0,\B)$ can be any number in $[0,R]$.

\end{rem}

We now show that the rate $\rho_*$ in Theorem \ref{mainthm} is optimal. The main point is to reduce the dynamics of $T$ near $\B$ to a map on $\Nor_{H_b}$.

\begin{prop}\label{propsharp}
For any $\eps\in(0,\rho_*)$ and $\delta>0$, there exists $H_0\in\M$ with $0<\dist(H_0,\B)<\delta$ such that
$$\dist\big(T^m(H_0),\B\big)\geq(\rho_*-\eps)^m\dist(H_0,\B)\qquad\textrm{for all } m\geq0.$$
Moreover, $H_0$ can be chosen so that $\dist\big(T^m(H_0),\B\big)\geq c\,\rho_*^m\dist(H_0,\B)$ for all $m\geq0$, where $c>0$ depends on $H_0$. In particular, the rate $\rho_*$ in Theorem \ref{mainthm} cannot be replaced by any $\lambda<\rho_*$.
\end{prop}

\begin{proof}
  We fix $H_b\in\B$ and work in a small neighbourhood of $H_b.$
    
   \medskip
\noindent {\bf STEP I:} 
\noindent 
In this step, we construct a map $\Psi:\Nor_{H_b}\to \Nor_{H_b}$ such that $$\dist\big(T^m(\Exp_{H_b}(y)),\B\big)=\norm{\Psi^m(y)},$$ as long as $\norm{\Psi(y)}, \dots \norm{\Psi^{m-1}(y)}$ are small. Define
$$F:\pp_{H_b}\times\Nor_{H_b}\To\M,\qquad F(u,y)=(e^{u/2})^*\Exp_{H_b}(y).$$
We prove that \begin{equation}\label{EqF}\pi(F(u,y))=\Exp_{H_b}(u) \,\,\,\, \textrm{and} \,\,\,\, \dist(F(u,y),\B)=\norm{y}.\end{equation}
Put $g=e^{u/2}\in G$. Since $g^*$ is an isometry of $\M$ preserving $\B$, we have $F(u,y)=\Exp_{g^*H_b}(g^{-1}yg)$ with $g^{-1}yg\in\Nor_{g^*H_b}$, and $g^*H_b=\Exp_{H_b}(u)$. Thus, Proposition \ref{propgeomB} implies \eqref{EqF}. By definition, the differential of $F$ at $(0,0)$ is $(u,y)\mapsto u+y$, which is an isomorphism onto $\Herm_{H_b}(W)=\pp_{H_b}\oplus\Nor_{H_b}$. Therefore the inverse function theorem implies that there exist a neighbourhood $U$ of $0$ in $\Herm_{H_b}(W)$ and smooth maps $u:U\to\pp_{H_b}$ and $v:U\to\Nor_{H_b}$ such that
$$\Exp_{H_b}(z)=F\big(u(z),v(z)\big)\qquad\textrm{for } z\in U.$$
Note that $v(z)=z$ for $z\in U\cap\Nor_{H_b}$ and $v(z)=0$ for $z\in U\cap\pp_{H_b}$. Hence the differential of $v$ at $0$ is the orthogonal projection onto $\Nor_{H_b}$.

Let $\Phi$ be as in the proof of Lemma \ref{lemlocal}, and choose $\delta_0>0$ such that $\Phi(y)\in U$ for $\norm{y}\leq\delta_0$. For $y\in\Nor_{H_b}$ with $\norm{y}\leq\delta_0$, define 
$$\Psi(y)=v(\Phi(y))\in\Nor_{H_b}.$$ By Theorem \ref{thmlin}(2), $\Psi$ is smooth, $\Psi(0)=0$, and its differential at $0$ is the restriction  of $d_{H_b}T$ to $\Nor_{H_b}$. Moreover, for $y \in \Nor_{H_b}$ close to the origin, we have
$$T\big(\Exp_{H_b}(y)\big)=\Exp_{H_b}(\Phi(y))=g^*\Exp_{H_b}(\Psi(y)),\qquad g=e^{u(\Phi(y))/2}\in G.$$
Since $T$ commutes with the action of $G$ and $\dist(\cdot,\B)$ is $G$-invariant, induction on $m$ gives
\begin{equation}\label{eqPsi}
\dist\big(T^m(\Exp_{H_b}(y)),\B\big)=\norm{\Psi^m(y)},
\end{equation}
as long as $\norm{\Psi^j(y)}\leq\delta_0$ for $j<m$.

 \medskip
\noindent {\bf STEP II:} 
\noindent
Put $A=d_{H_b}T|_{\Nor_{H_b}}$. Let $\Nor_{H_b}=E \oplus E^{\prime}$ where $E$ is the eigenspace of $A$ with eigenvalue $\rho_*$ and $E'$ its orthogonal complement in $\Nor_{H_b}$. We denote the orthogonal projections onto $E$ and $E'$ by $\Pi$ and $\Pi'$, and the largest eigenvalue of $A$ on $E'$ by $\rho'$. Let $\rho'=0$ if $E'=0$. Hence $\rho'<\rho_*$. Since $A$ is self-adjoint, we have $$\Pi A=\rho_*\Pi \,\,\,\,\, \textrm{and} \,\,\,\,\, \norm{\Pi'Ay}\leq\rho'\norm{\Pi'y}.$$ Choose $C_3>0$ such that $$\norm{\Psi(y)-Ay}\leq C_3\norm{y}^2, \,\, \textrm{for}\,\, \norm{y}\leq\delta_0.$$ Choose    $\delta_1=\min(\delta_0,\delta, \frac{\eps}{3C_3},\frac{\rho_*-\rho'}{3C_3},\frac{1-\rho_*}{3C_3})$  and let
$$S=\set{y\in\Nor_{H_b}\mid \norm{y}\leq\delta_1,\ \norm{\Pi'y}\leq\norm{\Pi y}}.$$
For $y\in S$, we have $\norm{y}\leq\sqrt2\norm{\Pi y}$, and therefore $C_3\norm{y}^2\leq\sqrt2C_3\delta_1\norm{\Pi y}$. Hence
\begin{align*}
\norm{\Psi(y)}&\leq(\rho_*+C_3\delta_1)\norm{y}\leq\norm{y},\\
\norm{\Pi\Psi(y)}&\geq(\rho_*-\sqrt2C_3\delta_1)\norm{\Pi y}\geq(\rho_*-\eps)\norm{\Pi y},\\
\norm{\Pi'\Psi(y)}&\leq\rho'\norm{\Pi'y}+\sqrt2C_3\delta_1\norm{\Pi y}\\
&\leq(\rho_*-\sqrt2C_3\delta_1)\norm{\Pi y}\leq\norm{\Pi\Psi(y)}.
\end{align*}
Therefore $\Psi(S)\subseteq S$. Since $\B\neq\M$, we have $E\neq0$. Choose $y_0\in E$ with $0<\norm{y_0}<\delta_1$ and put $H_0=\Exp_{H_b}(y_0)$. Then $\dist(H_0,\B)=\norm{y_0}<\delta$ and $\Psi^m(y_0)\in S$ for all $m$. By \eqref{eqPsi}, we get
$$\dist\big(T^m(H_0),\B\big)=\norm{\Psi^m(y_0)}\geq\norm{\Pi\Psi^m(y_0)}\geq(\rho_*-\eps)^m\norm{y_0}.$$
Moreover, put $y_m=\Psi^m(y_0)\in S$. The above estimates give $\norm{y_m}\leq(\rho_*+C_3\delta_1)^m\norm{y_0}$ and $\norm{\Pi y_{m+1}}\geq\rho_*\norm{\Pi y_m}-C_3\norm{y_m}^2\geq\big(\rho_*-\sqrt2C_3\norm{y_m}\big)\norm{\Pi y_m}$. Since $\rho_*+C_3\delta_1<1$ and $\sqrt2C_3\norm{y_m}<\eps<\rho_*$, we get
$$\dist\big(T^m(H_0),\B\big)\geq\norm{\Pi y_m}\geq c\,\rho_*^m\norm{y_0}.$$

Finally, let $\lambda<\rho_*$. Choose $\eps$ with $\rho_*-\eps>\lambda$ and $H_0$ as above. Since $H_\infty\in\B$, we have $\dist(T^m(H_0),H_\infty)\geq(\rho_*-\eps)^m\dist(H_0,\B)$, which is not bounded by $C\lambda^m$ for any constant $C$.
\end{proof}

\begin{rem}
The fact that $G$ acts on $\M$ by isometries plays a crucial role in the argument. It makes the constants in Lemmas \ref{lemloj}, \ref{lemlocal} and \ref{lementer} uniform along the noncompact submanifold $\B$, and it reduces the dynamics near $\B$ to the map $\Psi$ in the proof of Proposition \ref{propsharp}. Note that the drift of the sequence along $\B$ is then controlled by the distance to $\B$ (Lemma \ref{lemlocal}(2)), which decays exponentially.
\end{rem}

\section{Examples}\label{secexample}

In this section, we compute $\rho_*$ for $\CP^n$ and $\CP^1\times\CP^1$. As in the Remark \ref{remberezin}, we have   $d_HT=Q\circ(1-\Delta_h)\circ\sigma$.

\begin{lem}\label{lemeigen}
Let $H\in\B$, $h=\FS(H)$ and $a\in\Herm_H(W)$. Suppose that $\mathfrak{B}(\hat a)=\beta \hat a$ and $-\Delta_h \hat a=\gamma \hat a$. Then $d_HT(a)=\beta(1+\gamma)a$. 
\end{lem}

\begin{proof}
Let $\sigma(a)=\hat a.$ We have $\sigma\big(Q(\hat a)\big)=\mathfrak{B}\hat a=\beta\sigma(a)=\sigma(\beta a)$. By the proof of Theorem \ref{thmlin}(4), $\sigma$ is injective. Hence $Q(\sigma(a))=\beta a$. Therefore
$$d_HT(a)=Q\big((1-\Delta_h)\sigma(a)\big)=(1+\gamma)\,Q(\sigma(a))=\beta(1+\gamma)a.$$
\end{proof}

\subsection{Projective spaces}\label{secproj}

Let $X=\CP^n$ and $L=\mathcal{O}(d)$ with $d\geq2$. Note that for $d=1$ we have $\B=\M$. In this case $W=H^0(\CP^n,\mathcal{O}(d))\cong\Sym^d\C^{n+1}$, and $G$ is the image of $GL(n+1,\C)$. Let $H_b$ be the inner product on $W$ induced by the standard inner product on $\C^{n+1}=H^0(\CP^n,\mathcal{O}(1))$. Since $H_b$ is $U(n+1)$-invariant, so is $T(H_b)$ by Lemma \ref{lemequiv}. Since $W$ is an irreducible representation of $U(n+1)$, Schur's lemma gives $T(H_b)=cH_b$ for some $c>0$, and $c=1$ by Lemma \ref{lemtrace}. Hence $H_b$ is balanced. Therefore, by Proposition \ref{propB}, an inner product on $W$ is balanced if and only if it is induced by a Hermitian inner product on $\C^{n+1}$. Put $h=\FS(H_b)$. Then $\omega_{h}=d\,\omega_{FS}$, where $\omega_{FS}=\ddbar\log|x|^2$ is the Fubini--Study form on $\CP^n$. Let $\Delta$ be the Laplacian of $\omega_{FS}$, normalized as in Theorem \ref{thmlin}, so that $\Delta_h=d^{-1}\Delta$. The eigenvalues of $-\Delta$ are $\ell(\ell+n)$, $\ell\geq0$, and we denote by $E_\ell$ the corresponding real eigenspaces. The $E_\ell\otimes\C$ are pairwise inequivalent irreducible representations of $SU(n+1)$, and the functions $f(x,\bar x)/|x|^{2d}$ with $f$ bihomogeneous of bidegree $(d,d)$ form the space $\bigoplus_{\ell=0}^d E_\ell\otimes\C$ (c.f.~\cite{BGM}). Its complex dimension is $N^2=\dim_\C\mathrm{End}(W)$.

Let $s_1,\dots,s_N$ be an $H_b$-orthonormal basis of $W$. Since $H_b$ is $U(n+1)$-invariant, $\sum|s_i(x)|^2$ is a constant multiple of $|x|^{2d}$. For $a\in\Herm_{H_b}(W)$, we have
$$\sigma(a)([x])=\frac{\sum_{i,j}a_{ij}\,s_i(x)\overline{s_j(x)}}{\sum_k|s_k(x)|^2},$$
where $a_{ij}=H_b(as_j,s_i)$. Hence $|x|^{2d}\sigma(a)$ is bihomogeneous of bidegree $(d,d)$. Since $\sigma$ is injective and commutes with the action of $U(n+1)$, it is therefore an equivariant isomorphism from $\Herm_{H_b}(W)$ onto $\bigoplus_{\ell=0}^dE_\ell$. For $0\leq \ell\leq d$, let $U_{\ell}=\sigma^{-1}(E_\ell)$, i.e. the subspace of those $a\in\Herm_{H_{b}}(W)$ for which $\hat a$ is an eigenfunction of $-\Delta$ with eigenvalue $\ell(\ell+n)$. Then
$$\Herm_{H_{b}}(W)=\bigoplus_{\ell=0}^{d}U_{\ell},$$
and the $U_{\ell}\otimes\C$ are pairwise inequivalent irreducible representations of $SU(n+1)$.

The following lemma is contained in \cite[Lemma 3.7 and its proof]{Se2}, and we refer the reader to it for details. Note that, with respect to $\frac{N}{V}d\mu_h$, the integral kernel of $\mathfrak{B}$ is $|\langle x,y\rangle|^{2d}/(|x|^{2d}|y|^{2d})$, which only depends on the distance between $[x]$ and $[y]$. Hence $\mathfrak{B}$ commutes with $U(n+1)$, and it acts on each $E_\ell$ by a scalar.

\begin{lem}\label{lembeta}
The Berezin transform of $h=\FS(H_b)$ acts on $E_\ell$ by the scalar
$$\beta_\ell=\frac{d!\,(d+n)!}{(d-\ell)!\,(d+n+\ell)!}.$$
\end{lem}

\begin{thm}\label{thmspec}
For $0\leq\ell\leq d$, the map $d_{H_b}T$ acts on $U_\ell$ by the scalar
\begin{equation}\label{eqck}
\rho_{\ell}=\frac{d!\,(d+n)!}{(d-\ell)!\,(d+n+\ell)!}\Big(1+\frac{\ell(\ell+n)}{d}\Big).
\end{equation}
Moreover, $\rho_{0}=\rho_{1}=1$, $U_{0}\oplus U_{1}=T_{H_{b}}\B$, and $\rho_{\ell}<\rho_{2}$ for $\ell\geq3$. In particular,
$$\rho_{*}=\rho_{2}=1-\frac{(n+2)(n+3)}{(d+n+1)(d+n+2)}.$$
\end{thm}

\begin{proof}
Let $a\in U_\ell$. Then $\sigma(a)\in E_\ell$ is an eigenfunction of $-\Delta_h=-d^{-1}\Delta$ with eigenvalue $\ell(\ell+n)/d$, and of $\mathfrak{B}$ with eigenvalue $\beta_\ell$ by Lemma \ref{lembeta}. Therefore Lemma \ref{lemeigen} gives \eqref{eqck}. Simple calculation gives $\rho_0=\rho_1=1>\rho_2>\dots>\rho_d$. Note that the eigenspace of $d_HT $ associated to eigenvalue $1$ is $T_{H_b}\B=U_0\oplus U_1$. Finally,
\[
\rho_*=\rho_2
=1-\frac{(n+2)(n+3)}{(d+n+1)(d+n+2)}.\]

\end{proof}

\begin{rem}\label{remsherman}
For $n=1$, Theorem \ref{thmspec} gives $\rho_*=\frac{(d-1)(d+6)}{(d+2)(d+3)}$. This is the rate of convergence obtained numerically by Sherman for toric metrics on $(\CP^1,\mathcal{O}(d))$ (\cite[(9)]{Sh}). For general $n$, we have $1-\rho_*\sim(n+2)(n+3)d^{-2}$ as $d\to\infty$. Hence the convergence of the algorithm slows down for large powers of the polarization. Note that Theorem \ref{mainthm} does not require the initial metric to be $S^1$-invariant.
\end{rem}

\begin{rem}\label{remthompson}
Let $n=1$ and $d=2$, and let $s_0,s_1,s_2$ be the $H_b$-orthonormal basis of $W$ obtained by rescaling $z_0^2$, $z_0z_1$, $z_1^2$. In this basis, $T_{H_b}\B=U_0\oplus U_1$ is spanned by $I$, $J=\mathrm{diag}(1,0,-1)$ and two Hermitian matrices with vanishing diagonal, namely the images of the Hermitian $2\times2$ matrices under the derived action. By Theorem \ref{thmspec}, $\rho_*=\frac25$ and $d_{H_b}T=\frac25I+\frac35\Pi$, where $\Pi$ is the orthogonal projection onto $T_{H_b}\B$. Let $e=\mathrm{diag}(1,0,0)$. Then $\Pi(e)=\frac13I+\frac12J$, and therefore
$$d_{H_b}T(e)=\mathrm{diag}\Big(\frac{9}{10},\frac15,-\frac{1}{10}\Big).$$
Since $e\geq0$ and $d_{H_b}T(e)$ has a negative eigenvalue, $T$ is not order preserving: for small $\eps>0$, we have $\Exp_{H_b}(\eps e)\geq H_b$, while $T\big(\Exp_{H_b}(\eps e)\big)\not\geq H_b=T(H_b)$. Moreover, for $Y=I-2e$ we have $d_{H_b}T(Y)=\mathrm{diag}\big(-\frac45,\frac35,\frac65\big)$. Hence, for $H_\eps=\Exp_{H_b}(\eps Y)$, we get
$$d_{Th}(H_\eps,H_b)=\eps,\qquad d_{Th}\big(T(H_\eps),T(H_b)\big)=\frac65\,\eps+O(\eps^2).$$
Therefore $T$ is not distance non-increasing with respect to the Thompson metric, even near a balanced metric. On the other hand, the operator norm of $d_{H_b}T$ with respect to the Riemannian metric is $1$ by Theorem \ref{thmlin}. Note also that, far from $\B$, the map $T$ can increase the Riemannian distance to its limit (\cite[Section 4]{Sh}).
\end{rem}

\subsection{The product $\CP^1\times\CP^1$}

Let $X=\CP^1\times\CP^1$ and $L=\mathcal{O}(d,d)=\pi_1^*\mathcal{O}(d)\otimes\pi_2^*\mathcal{O}(d)$ with $d\geq1$. Then $W=W_1\otimes W_1$, where $W_1=H^0(\CP^1,\mathcal{O}(d))$. Let $H_1$ be the balanced inner product of Section \ref{secproj} on $W_1$ (with $n=1$), $h_1=\FS(H_1)$, and $H_b=H_1\otimes H_1$. If $s_1,\dots,s_{d+1}$ is an $H_1$-orthonormal basis of $W_1$, then $s_i\otimes s_j$ is an $H_b$-orthonormal basis of $W$. Hence $h=\FS(H_b)=h_1\boxtimes h_1$, $\omega_h=\pi_1^*\omega_{h_1}+\pi_2^*\omega_{h_1}$ and $d\mu_h=d\mu_{h_1}\times d\mu_{h_1}$. Therefore $\Hilb(h)=\Hilb(h_1)\otimes\Hilb(h_1)$, and $H_b$ is balanced. Moreover, with the notation of Remark \ref{remberezin}, we have
$$Q(f_1\otimes f_2)=Q_1(f_1)\otimes Q_1(f_2),\qquad \sigma(a_1\otimes a_2)=\sigma_1(a_1)\otimes\sigma_1(a_2),$$
$$\Delta_h=\Delta_{h_1}\otimes1+1\otimes\Delta_{h_1},$$
where $Q_1$ and $\sigma_1$ are the corresponding operators for $(\CP^1,\mathcal{O}(d))$ at $H_1$. In particular, the Berezin transform of $h$ is $\mathfrak{B}_1\otimes\mathfrak{B}_1$. Let $U_\ell\subset\Herm_{H_1}(W_1)$ be as in Section \ref{secproj}. Since $\Herm_{H_b}(W)=\Herm_{H_1}(W_1)\otimes_\R\Herm_{H_1}(W_1)$, we have
$$\Herm_{H_b}(W)=\bigoplus_{\ell,k=0}^dU_\ell\otimes U_k.$$

\begin{prop}\label{propproduct}
For $0\leq\ell,k\leq d$, the map $d_{H_b}T$ acts on $U_\ell\otimes U_k$ by the scalar
$$\rho_{\ell,k}=\beta_\ell\beta_k\Big(1+\frac{\ell(\ell+1)+k(k+1)}{d}\Big),\qquad \beta_\ell=\frac{d!\,(d+1)!}{(d-\ell)!\,(d+\ell+1)!}.$$
Moreover, $T_{H_b}\B=U_0\otimes U_0\oplus U_1\otimes U_0\oplus U_0\otimes U_1$, and
$$\rho_*=\rho_{1,1}=\frac{d(d+4)}{(d+2)^2}=1-\frac{4}{(d+2)^2}.$$
\end{prop}

\begin{proof}
Let $a_1\in U_\ell$ and $a_2\in U_k$. By Lemma \ref{lembeta}, $\sigma(a_1\otimes a_2)=\sigma_1(a_1)\otimes\sigma_1(a_2)$ is an eigenfunction of $\mathfrak{B}_1\otimes\mathfrak{B}_1$ with eigenvalue $\beta_\ell\beta_k$, and of $-\Delta_h$ with eigenvalue $(\ell(\ell+1)+k(k+1))/d$. Hence the first statement follows from Lemma \ref{lemeigen}.

 Let $\rho_\ell$ be the numbers given  by Theorem \ref{thmspec} for $n=1$. It is easy to show that 
$$\frac{\rho_{\ell+1,k}}{\rho_{\ell,k}}\leq \frac{\rho_{\ell+1}}{\rho_\ell},$$
with equality if and only if $k=0$. Therefore,
$\rho_{0,0}=\rho_{1,0}=\rho_{0,1}=1$ which implies that  $T_{H_b}\B$ is the eigenspace of $d_{H_b}T$ with eigenvalue $1$, 
\end{proof}

Note that the slowest direction $U_1\otimes U_1$ is a mixed mode, which has no counterpart on either factor. Hence the convergence on $\CP^1\times\CP^1$ is slower than on $\CP^1$. For instance, for $d=2$ we get $\rho_*=\frac34$, while $\rho_*=\frac25$ for $(\CP^1,\mathcal{O}(2))$.

\subsection{Large powers of the polarization}\label{secasym}

Let $\omega_\infty\in2\pi c_1(L)$ be a constant scalar curvature K\"ahler metric. Suppose that $L^k$ admits balanced inner products $H_k$ for $k\geq k_0$ with $k^{-1}\omega_{\FS(H_k)}\to\omega_\infty$ in $C^\infty$. We denote by $\rho_*(k)$ the number $\rho_*$ of Theorem \ref{thmlin}(5) for $(X,L^k)$, i.e. the rate in Theorem \ref{mainthm}. Let $\mu_1$ be the smallest positive eigenvalue of the Lichnerowicz operator
$$\mathcal{D}^*\mathcal{D}f=\Delta^2f+R^{j\bar k}f_{j\bar k}$$
of $\omega_\infty$, where $\Delta$ is the Laplacian of $\omega_\infty$ normalized as in Theorem \ref{thmlin} and $R^{j\bar k}$ is the Ricci tensor of $\omega_\infty$ with indices raised. One expects that
\begin{equation}\label{eqrate}
1-\rho_*(k)=\frac{\mu_1}{2k^2}+O(k^{-3}).
\end{equation}
When $\Aut(X,L)/\C^*$ is discrete, Ioos and Polterovich prove \eqref{eqrate} with $o(k^{-2})$ in place of $O(k^{-3})$ (\cite[Theorem 1.2(2)]{IP}). Note that in \cite{IP} the K\"ahler form is normalized by $\omega=\frac{\omega_\infty}{2\pi}\in c_1(L)$, and the leading term is written as $\frac{\mu_1(D)}{8\pi k^2}$, where $D$ is the operator of \cite[(1.18)]{IP}. By \cite[(4.43)]{IP}, we have $D=\pi^{-1}\mathcal{D}^*\mathcal{D}_\omega=4\pi\mathcal{D}^*\mathcal{D}$, where $\mathcal{D}^*\mathcal{D}_\omega$ is the operator above computed with respect to $\omega$. Hence $\mu_1(D)=4\pi\mu_1$, and the two normalizations agree.

The examples above satisfy \eqref{eqrate}. For $(\CP^n,\mathcal{O}(1))$ with $\omega_\infty=\omega_{FS}$, we have $\mathrm{Ric}(\omega_{FS})=(n+1)\omega_{FS}$. Hence $\mathcal{D}^*\mathcal{D}=\Delta(\Delta+n+1)$ acts on $E_\ell$ by $\ell(\ell+n)\big(\ell(\ell+n)-n-1\big)$, so that $\mu_1=2(n+2)(n+3)$. Theorem \ref{thmspec} gives
$$1-\rho_*(k)=\frac{(n+2)(n+3)}{k^2}+O(k^{-3})=\frac{\mu_1}{2k^2}+O(k^{-3}).$$
For $(\CP^1\times\CP^1,\mathcal{O}(1,1))$ with the product metric, $\mathcal{D}^*\mathcal{D}=\Delta(\Delta+2)$ acts on $E_\ell\otimes E_k$ by $\nu(\nu-2)$, where $\nu=\ell(\ell+1)+k(k+1)$. Hence $\mu_1=8$, attained on $E_1\otimes E_1$, and Proposition \ref{propproduct} gives $1-\rho_*(k)=4k^{-2}+O(k^{-3})=\mu_1/(2k^2)+O(k^{-3})$. Note that in both examples $\Aut(X,L)/\C^*$ is not discrete, and the slowest direction of $T$ corresponds to the eigenspace of $\mathcal{D}^*\mathcal{D}$ with eigenvalue $\mu_1$.

We do not know whether \eqref{eqrate} holds in general when $\Aut(X,L)/\C^*$ is not discrete. Moreover, the error term $O(k^{-3})$ is not known in general, even when $\Aut(X,L)/\C^*$ is discrete. Note also that if $\Aut(X,L)/\C^*$ is not discrete, the existence of balanced metrics for all $k\gg0$ is a genuine assumption and is not guaranteed by the existence of a cscK metric. In \cite{OSY}, Ono, Sano and Yotsutani construct a cscK polarization which is asymptotically Chow unstable and therefore does not admit balanced metrics for high powers (see Example \ref{exosy}).

\appendix
\section{Invariance of $\Zt$ and Futaki's obstruction}\label{appfutaki}

In Lemma \ref{lemequiv}, the invariance $\Zt(g^*H)=\Zt(H)$ for $g\in G$ is proved under the assumption that $(X,L)$ admits a balanced metric. In this appendix, we show that this assumption cannot be removed, and we explain the relation of this invariance with the obstruction to asymptotic Chow semistability of Mabuchi and Futaki (\cite{Ma}, \cite{Fu2}). Throughout this appendix, we do not assume that $(X,L)$ admits a balanced metric.

Let $\mathfrak{g}\subseteq\textrm{End}(W)$ be the Lie algebra of $G$. For $H\in\M$, we define the complex linear functional $F_H:\mathfrak{g}\to\C$ by
$$F_H(b)=\tr\big(bP(H)\big).$$
Let $h=\FS(H)$ and let $s_1,\dots,s_N$ be an $H$-orthonormal basis of $W$. By the definition of $\Hilb$, we have
$$F_H(b)=\frac{N}{V}\int_X\hat b\,d\mu_h-\tr(b),$$
where $\hat b(x)=\sum_{i,j}H(bs_j,s_i)\,h\big(s_i(x),s_j(x)\big)$. If $b$ is $H$-self-adjoint, then $\hat b=\mu_b\circ\iota_H$ is the function of Section \ref{subsecsecond}, where $\mu_b$ is as in Lemma \ref{lemsecond}.

\begin{prop}\label{propchar}
The following hold.
\begin{enumerate}
\item There exists a continuous homomorphism $c:G\to\R$ such that $\Zt(g^*H)=\Zt(H)+c(g)$ for all $g\in G$ and $H\in\M$.
\item The functional $F_H$ does not depend on $H\in\M$. We denote it by $F$. Moreover, for $b\in\mathfrak{g}$ and $t\in\R$, we have
$$c(e^{tb})=-\frac{2V}{N}\,t\,\textrm{Re}\,F(b).$$
\item The functional $F$ vanishes on the scalars and on $[\mathfrak{g},\mathfrak{g}]$, and $c$ vanishes on every compact subgroup of $G$.
\item The functional $\Zt$ is $G$-invariant if and only if $F=0$.
\item If $(X,L)$ admits a balanced metric, then $F=0$.
\end{enumerate}
\end{prop}

\begin{proof}
The first part of the proof of Lemma \ref{lemequiv} does not use the existence of a balanced metric, and it gives $(1)$. Note that $c$ is continuous, since $\Zt$ is continuous.

Let $b\in\mathfrak{g}$ and $H\in\M$, and put $H_t=(e^{tb})^*H$. Then
$$\frac{d}{dt} H_t(s,s')=H_t(bs,s')+H_t(s,bs')=H_t\big((b+b^{*_t})s,s'\big),$$
where $b^{*_t}$ is the adjoint of $b$ with respect to $H_t$. 
Since $P(H_t)$ is $H_t$-self-adjoint, we have $\tr\big(b^{*_t}P(H_t)\big)=\overline{\tr\big(bP(H_t)\big)}$. 
Therefore Lemma \ref{lemfirst} implies that
$$\frac{d}{dt}\Zt(H_t)=-\frac{V}{N}\tr\big((b+b^{*_t})P(H_t)\big)=-\frac{2V}{N}\,\textrm{Re}\,F_{H_t}(b).$$
On the other hand, $\Zt(H_t)=\Zt(H)+c(e^{tb})$ by $(1)$, and $t\mapsto c(e^{tb})$ is a continuous homomorphism from $\R$ to $\R$. Hence $c(e^{tb})=t\,\tilde{c}(b)$ for a real number $\tilde{c}(b)$, and we get $\textrm{Re}\,F_{H_t}(b)=-\frac{N}{2V}\tilde{c}(b)$ for all $t$. Therefore $\textrm{Re}\,F_H(b)$ does not depend on $H$. Since $F_H$ is complex linear and $\textrm{Im}\,F_H(b)=-\textrm{Re}\,F_H(\sqrt{-1}\,b)$, the same holds for $F_H(b)$. This proves $(2)$.

Since $\tr P(H)=0$, then  $F$ vanishes on the scalars. The map $b\mapsto \tilde{c}(b)$ is the differential at the identity of the Lie group homomorphism $c$. Hence it is a homomorphism of Lie algebras from $\mathfrak{g}$ to $\R$, and therefore it vanishes on $[\mathfrak{g},\mathfrak{g}]$. By $(2)$, $\textrm{Re}\,F$ vanishes on the complex subspace $[\mathfrak{g},\mathfrak{g}]$, and so does $F$. Finally, $c(K)$ is a compact subgroup of $\R$ for a compact subgroup $K\subseteq G$ and therefore, $c(K)=0.$  This proves $(3)$.

If $\Zt$ is $G$-invariant, then $c=0$, and $(2)$ gives $F=0$. Conversely, suppose that $F=0$. By $(2)$, $c$ vanishes on $e^{tb}$ for every $b\in\mathfrak{g}$, and therefore on the identity component $G_0$ of $G$. Since $G$ is the subgroup of $GL(W)$ preserving the subvariety $\iota(X)$ of $\mathbb{P}(W^*)$, it is an algebraic group, and therefore $G/G_0$ is finite. Hence for any $g\in G$ there exists $m\geq1$ with $g^m\in G_0$, and $m\,c(g)=c(g^m)=0$. This proves $(4)$.

Finally, if $H_b$ is balanced, then $P(H_b)=0$, and therefore $F=F_{H_b}=0$. This proves $(5)$.
\end{proof}

Proposition \ref{propchar} is closely related to the work of Mabuchi and Futaki recalled in Remark \ref{remchow}, and we include it for the convenience of the reader. Note that $(4)$ and $(5)$ give another proof of the last statement of Lemma \ref{lemequiv}. Moreover, $(2)$ and $(5)$ give an obstruction to the existence of balanced metrics. Indeed, if $F\neq0$, then $(X,L)$ admits no balanced metric, and $\Zt$ is not bounded from below on $\M$, since $c(e^{tb})$ is a nonconstant linear function of $t$ for some $b\in\mathfrak{g}$.

\begin{rem}\label{remchow}
Note that the functional $F$ is the derivative of the character by which $G$ acts on the Chow point of $X$. Indeed, let $g\in G\cap SL(W)$. Since $g$ preserves $\iota(X)$, we have $g\cdot\textrm{Chow}(\iota(X))=\chi(g)\,\textrm{Chow}(\iota(X))$ for a character $\chi:G\cap SL(W)\to\C^*$. The identity of Zhang recalled in Remark \ref{remKN} holds with any $H\in\M$ as the base point (\cite{Zh}, \cite{PS}). Indeed, choosing an $H$-orthonormal basis of $W$ identifies $(W,H)$ with $\C^N$ and its standard inner product. This identification maps $\iota(X)$ to a subvariety of $\CP^{N-1}$, and it maps the norm on the space of Chow forms induced by $H$ to the one induced by the standard inner product. Since the identity holds for every subvariety of $\CP^{N-1}$ with the standard inner product as the base point, it also holds with $H$ as the base point. Therefore $c(g)=\theta\log|\chi(g)|$, where $\theta>0$ depends only on $(X,L)$. Since $\Zt$ is invariant under scaling, we get that $F=0$ if and only if $|\chi|=1$ on $G\cap SL(W)$. Suppose that $F\neq0$. By Proposition \ref{propchar}(3), there exists $b\in\mathfrak{g}$ with $\tr(b)=0$ and $\textrm{Re}\,F(b)\neq0$. Then $|\chi(e^{tb})|$ tends to $0$ as $t\to\infty$ or as $t\to-\infty$, so that $0$ lies in the closure of the $SL(W)$-orbit of the Chow point. Hence, if the Chow point of $\iota(X)$ is semistable, then $F=0$. 

In \cite{Ma}, Mabuchi shows that the vanishing of these characters for the polarizations $L^k$, $k\gg0$, is necessary for the asymptotic Chow semistability of $(X,L)$. In \cite{Fu2}, Futaki expresses this obstruction in terms of a family of integral invariants, the first of which is the classical Futaki invariant of \cite{Fu1}. 
\end{rem}

We now compute $F$ for toric manifolds. Let $(X,L)$ be a polarized toric manifold of dimension $n$, and let $P\subseteq\R^n$ be the polytope associated with $L$. Then $W$ has a basis of sections $s_m$, $m\in P\cap\mathbb{Z}^n$, on which the torus $(\C^*)^n\subseteq\Aut(X,L)$ acts by the characters $m$. In particular, $N=\#(P\cap\mathbb{Z}^n)$. For $\xi\in\R^n$, we denote by $b_\xi\in\mathfrak{g}$ the endomorphism defined by $b_\xi s_m=\langle\xi,m\rangle s_m$, and by $\textrm{bar}(P)$ the barycentre of $P$ with respect to the Lebesgue measure.

\begin{prop}\label{proptoric}
For any $\xi\in\R^n$, we have
$$F(b_\xi)=N\Big\langle\xi,\ \textrm{bar}(P)-\frac1N\sum_{m\in P\cap\mathbb{Z}^n}m\Big\rangle.$$
In particular, $F$ vanishes on the Lie algebra of $(\C^*)^n$ if and only if the average of the lattice points of $P$ is the barycentre of $P$.
\end{prop}

\begin{proof}
By Proposition \ref{propchar}(2), we may compute $F$ with the inner product $H$ for which the sections $s_m$ are orthonormal. Put $h=\FS(H)$. Since $\sum_m|s_m|_h^2=1$, the map
$$\mu=\sum_m m\,|s_m|^2_h:X\To\R^n$$
is, up to a constant factor, a moment map for the action of the compact torus $T^n$ on $(X,\omega_h)$, and its image is the convex hull $P$ of the weights $m$. Moreover, since $b_\xi$ is diagonal in the basis $s_m$, we have
$$F(b_\xi)=\sum_m\langle\xi,m\rangle\big(T(H)(s_m,s_m)-1\big)=\frac{N}{V}\int_X\langle\xi,\mu\rangle\,d\mu_h-\sum_m\langle\xi,m\rangle.$$
By the theorem of Duistermaat and Heckman (\cite{DH}), the push-forward of $d\mu_h$ by $\mu$ is a constant multiple of the Lebesgue measure on $P$. Hence $\frac1V\int_X\langle\xi,\mu\rangle\,d\mu_h=\langle\xi,\textrm{bar}(P)\rangle$, which proves the formula.
\end{proof}

In \cite[Theorem 1.2]{On}, Ono shows that if $(X,L)$ is Chow semistable, then the average of the lattice points of $P$ is the barycentre of $P$. Proposition \ref{proptoric} and Remark \ref{remchow} give another proof of this fact, which relies on the identity of Zhang recalled in Remark \ref{remchow}.

\begin{ex}\label{exhom}
Let $X=\CP^n$ and $L=\mathcal{O}(d)$. Then $\mathfrak{g}$ is spanned by the scalars and the image of $\mathfrak{sl}(n+1,\C)$. Since $\mathfrak{sl}(n+1,\C)=[\mathfrak{sl}(n+1,\C),\mathfrak{sl}(n+1,\C)]$, Proposition \ref{propchar}(3) gives $F=0$. Similarly, $F=0$ for $(\CP^1\times\CP^1,\mathcal{O}(a,b))$ with $a,b\geq1$. In the toric picture, $P$ is $d$ times the standard simplex, respectively a rectangle, and by symmetry the average of its lattice points is its barycentre. Note that balanced metrics exist in these cases (c.f. Section \ref{secexample}).
\end{ex}

\begin{ex}\label{exblowup}
Let $X$ be the blow-up of $\CP^2$ at a point, and let $L=-kK_X$ with $k\geq1$. Then $X$ is toric, and the polytope of $L$ is $kP$, where
$$P=\set{x\in\R^2\mid x_1\geq-1,\ x_2\geq-1,\ -1\leq x_1+x_2\leq1}$$
is the quadrilateral with vertices $(-1,0)$, $(-1,2)$, $(2,-1)$ and $(0,-1)$. A direct computation gives
\begin{align*}
N&=(2k+1)^2,\qquad \textrm{bar}(kP)=\frac{k}{12}(1,1),\\
\sum_{m\in kP\cap\mathbb{Z}^2}m&=\frac{k(k+1)(2k+1)}{6}(1,1).
\end{align*}
Hence, by Proposition \ref{proptoric},
$$F(b_\xi)=-\frac{k(2k+1)}{12}\,(\xi_1+\xi_2).$$
Therefore, for every $k\geq1$, the functional $\Zt$ on the space $\M$ associated with $(X,-kK_X)$ is neither $G$-invariant nor bounded from below, and $(X,-kK_X)$ admits no balanced metric. By Remark \ref{remchow}, the Chow point of $(X,-kK_X)$ is unstable for every $k\geq1$. Note that the coefficient of $k^2$ in $F(b_\xi)$ is
$$-\frac{\xi_1+\xi_2}{6}=\frac12\Big\langle\xi,\ |\partial P|\,\textrm{bar}(P)-\int_{\partial P}x\,d\sigma\Big\rangle,$$
where $d\sigma$ is the measure on $\partial P$ normalized so that $|\partial P|=8$ and $\int_{\partial P}x\,d\sigma=(1,1)$. By \cite{D4}, this is a multiple of the classical Futaki invariant of $(X,-K_X)$, which does not vanish by \cite{Fu1}.
\end{ex}

\begin{ex}\label{exosy}
In \cite{OSY}, Ono, Sano and Yotsutani construct a seven-dimensional toric Fano manifold $X$ which admits a K\"ahler--Einstein metric, such that $(X,-K_X)$ is asymptotically Chow unstable. In their example, the obstruction of Mabuchi and Futaki recalled in Remark \ref{remchow} does not vanish (\cite[Theorem 1.6]{OSY}). Therefore, by Proposition \ref{propchar} and Remark \ref{remchow}, for $k\gg0$ the functional $\Zt$ on the space $\M$ associated with $(X,-kK_X)$ is not $G$-invariant, and $(X,-kK_X)$ admits no balanced metric, although $X$ admits a K\"ahler--Einstein metric. In particular, the assumption in Section \ref{secasym} that balanced metrics exist for $k\gg0$ is not implied by the existence of a constant scalar curvature K\"ahler metric when $\Aut(X,L)/\C^*$ is not discrete.
\end{ex}

\end{document}